\documentclass[11pt]{article}
\usepackage{a4}
\usepackage{graphicx}
\usepackage{parskip}

\usepackage{amsfonts}
\usepackage{natbib}
\usepackage{algorithm}
\usepackage{algorithmic}

\usepackage{latexsym}

\newtheorem{theorem}{Theorem}[section]

\newtheorem{lemma}{Lemma}[section]
\newtheorem{corollary}{Corollary}[section]
\newtheorem{remark}{Remark}[section]

\newtheorem{definition}{Definition}[section]

\newcommand{\R}{\mathbb{R}}
\newcommand{\Nat}{\mathbb{N}}

\newcommand{\argmin}{\mathop{\arg \min}\limits}

\newcommand{\PP} {{  \rm I\hskip-0.22em P}}

\newcommand{\EE} {{\rm I\hskip-0.48em E}}

\usepackage{a4}

\usepackage{graphicx}

\usepackage{parskip}

\usepackage{amsfonts}

\usepackage{natbib}

\usepackage{algorithm}

\usepackage{algorithmic}

\usepackage{amsmath, amssymb}

\usepackage{latexsym}

\begin{document}

\centerline{\Large Some exercises with the Lasso and its compatibility constant}

\vskip .1in
\centerline{Sara van de Geer}

\centerline{January 12, 2017}

{\bf Abstract} We consider the Lasso for a noiseless experiment where one
has observations $X \beta^0$ and uses the penalized version of basis pursuit.
We compute for some special designs the compatibility constant,
a quantity closely related to the restricted eigenvalue. We moreover
show the dependence of the (penalized) prediction error on this compatibility constant.
This exercise illustrates that compatibility is necessarily entering into
the bounds for the (penalized) prediction error and that the bounds in the literature therefore
are - up to constants - tight. We also give conditions that show that in the noisy
case the dominating term for the prediction error is given by the
prediction error of the noiseless case.



\section{Introduction}
Let $X \in \R^{n \times p}$ be an $n \times p$ matrix and $\beta^0 \in \R^p$ be a fixed vector.
We consider the Lasso for the noiseless case
$$ \beta^* := \arg \min_{\beta \in \R^p } {\cal L}(\beta) , $$
with\footnote{In the noiseless case the results apply when $\| X \beta \|_2^2 $ ($\beta \in \R^p$) is replaced
by any other quadratic form $\beta^T \Sigma \beta$ ($\beta \in \R^p$) with $\Sigma$ a given
$p \times p$ matrix. The ``sample size" $n$ is playing the role of the rank of $\Sigma$.} 
$${\cal L} (\beta) :=\| X (\beta - \beta^0) \|_2^2 + 2 \lambda \| \beta \|_1 . $$
Aim in this note is to show that the upper bounds for $\| X ( \beta^* - \beta^0) \|_2^2 $
given in the literature (see Section \ref{upper-bound.section} for some references)
are also lower bounds, in the sense that there are designs where an upper bound is
tight, possibly up to constants. The upper bounds that we consider
depend on the so-called compatibility constant $\hat \phi^2  (S)$ which we define
in Definition \ref{compatibility.definition} below. 
In \cite{zhang2014lower} it is shown that for a given sparsity level, there is a design and a lower bound
for the mean prediction error  in the noisy case, that holds for any polynomial time algorithm. This lower bound is close to the known upper bounds and
in particular shows that compatibility conditions or restricted eigenvalue conditions cannot be avoided.
Our aim is  to make this visible for the Lasso by presenting some explicit expressions. 
This helps to understand why compatibility is playing a crucial role and also to understand the concept
itself. Our results follow from straightforward computation for some special
cases of design. 

We will show that the upper bounds involving
compatibility constants (here given in Section \ref{upper-bound.section}) match the lower
bounds  ``up to constants" or even ``asymptotically exactly" for certain designs. The designs we consider are in our view not 
atypical. Therefore, our conclusion is that there is not much space for improvement of
the existing upper bounds. 

Note that we consider a noiseless version of the Lasso. When examining lower bounds this is
reasonable, as one may expect that adding noise will not improve the performance of
the Lasso. We will moreover show in Section \ref{noisy.section} that for certain designs,
the ``bias" $\| X( \beta^*- \beta^0) \|_2^2$ of the noisy Lasso  is the dominating term,
so that bounds for the noiseless case immediately carry over to the noisy
case.

In order to be able to define the compatibility constant $\hat \phi^2 (S)$ we introduce here some notation.
For $S \subset \{1 , \ldots , p \}$  and a vector $\beta \in \R^p$ let
$ \beta_{j,S} := \beta_j {\rm l} \{ j \in S \} \in \R^p$. We apply the same notation for the $|S|$-dimensional vector 
$\{ \beta_j \}_{j \in S}$. We moreover write $\beta_{-S} := \beta_{S^c}$ where $S^c$ is the the complement
of the set $S$. If $S$ consists of a single variable, say $S= \{j \}$ we write 
$\beta_{-S} =: \beta_{-j} $.

\begin{definition}\label{compatibility.definition} 
The compatibility constant (see \cite{vandeGeer:07a} or  \cite{vdG2016} and its references) is
$$ \hat \phi^2 (L,S) := \min \biggl \{ |S| \| X \beta \|_2^2 : \ \| \beta_S \|_1 = 1 , \ \| \beta_{-S} \|_1 \le L \biggr \} .$$
The constant $L \ge 1$ is called a stretching factor. For $L=1$ we write $\hat \phi^2 (S):= \hat \phi^2 (1, S)$. 
When $S = \{ 1 , \ldots , p \}$ we let $\hat \phi^2 (S) := \min \{ |S| \| X \beta_S \|_2^2   :\ 
\| \beta_S \|_1 = 1 \} $.  For $S=\emptyset$ we set $|S|/ \hat \phi^2 ( S )= 0 $. 
\end{definition} 

The compatibility constant $\hat \phi^2 (L,S)$ with stretching constant $L>1$ can play a role when
considering the noisy situation. In this paper however, we mainly study the noiseless case
and take $L=1$. A noisy case where $L$ can be taken equal to 1 is considered in Section \ref{noisy.section}. 

It is sometimes helpful to consider $\hat \Gamma^2 (S) := |S| / \hat \phi^2 (S) $ as the 
{\it effective sparsity}\footnote{A better terminology is perhaps to call
$\hat \Gamma^2 (S)$ the effective {\it non}-sparsity} at the set $S$ (\cite{vdG2016}).  Two sets 
should be compared in terms of their effective sparsity rather than in terms of
their compatibility constants, in the sense that that we prefer sets $S$ with
$ \hat \Gamma^2 (S) $ small.

The compatibility constant $\hat \phi^2 (S)$
depends on the set $S$ and clearly also on the design $X$
through the Gram matrix $\hat \Sigma:= X^T X $. We express the latter dependence 
in our notation by the ``hat".  This is a habit coming from the case of random design,
where $\hat \Sigma$ is an estimator of $\EE \hat \Sigma$ (in statistics, estimators are
commonly denoted with a ``hat"). However, to avoid a cumbersome notation, not {\it all}
quantities depending on $X$ with be furnished with a ``hat". 

\subsection{Notation}\label{notation.section}

Let $X_j$ denote the $j$-th column of $X$ ($j=1 , \ldots , p $). The Gram matrix is
$\hat \Sigma := X^T X $. 

The active set (or support set) of $\beta^0$ is $S_0 := \{ j: \ \beta_j^0 \not= 0 \}$.
If $j \in S_0$ we call $j$ - or $X_j$ - an active variable. 
When $\hat \phi^2 (S_0)>0$ one says that the null space property holds  
(\cite{donoho2005neighborliness}). 
The cardinality of $S_0$ is denoted by $s_0 := | S_0 |$. 
We moreover write the cardinality of the set $S_0^c$ of inactive variables as $m_0 := p-s_0$.

\subsection{Organization of the paper}

Section \ref{noisy.section} shows how 
the results for the noiseless case carry over to the
noisy case when the Gram matrix (or an approximation thereof)
has bounded maximal eigenvalue and $\sqrt n \lambda$ is large
($\sqrt n \lambda \rightarrow \infty$). Such a choice for the tuning parameter
$\lambda$ corresponds to $p$ large, as can be the case in most of the following
sections (Sections \ref{p=s0=2N.section}, \ref{trivial.section}, \ref{s0=2N-m0=1.section},
\ref{s0=2-m0.section},
 \ref{s0=m0=2N.section}, and the last result of Section \ref{Furthers0=2.section}).
Section \ref{upper-bound.section} states some upper bounds for
the (penalized) prediction error of the noiseless Lasso. 
 These bounds are not novel, but as constants may now come into
play, we have re-derived them with an eye on the constants for the special situation with no noise.
Section \ref{design.section} has some considerations about the design:
we assume it to be ``fair" as defined there. Then, in the rest of the paper, we take the first
two variables as being among the active ones. In Section \ref{first2.section} we present the structure
(design and coefficients) for these first two variables. 
Section \ref{p=s0=2.section} considers the case $p=s_0=2$: it has no inactive variables.
This is extended in 
Section \ref{p=s0=2N.section} where $p=s_0=2N$ (for some $N \in \Nat$) is even. 
The next step is to start adding inactive variables. 
Section \ref{trivial.section} contains a trivial case,
where the inactive variables are orthogonal to the active ones.
Section \ref{s0=2-m0=1.section} has $s_0=2$ and $m_0=1$ and the single inactive variable
is a linear combination of the two active ones plus an orthogonal term: the active variables
are so to speak the ``parents" of the inactive one.
 Section \ref{s0=2N-m0=1.section} extends this to $s_0=2N$ even and $m_0=1$.
 Section \ref{s0=2-m0.section} returns to the case $s_0=2$, but now $m_0$ is arbitrary. The
 active variables are again ``parents" of all the inactive ones.
In Section \ref{s0=2=m0.section} we take $s_0$ as well as $m_0$ equal to 2,
but now part of the correlation between the two inactive variables 
is unique to those two, i.e., their correlation is not solely due to having the active ones as 
common ``parents".
Section \ref{s0=m0=2N.section} extends this to $s_0=m_0=2N$.
In Section \ref{Furthers0=2.section} the active variables are a linear combination of the inactive ones
plus orthogonal term: the inactive ones are now presented as the ``parents" of the active ones
instead of the other way around. 
Section \ref{proofs.section} contains the proofs.

For a symmetric matrix $A$ we let $\Lambda_{\rm min} (A)$ be its smallest and 
$\Lambda_{\rm max} (A)$ be its largest eigenvalue.

For two constants $u$ and
$v$ we let $u \vee v := \max\{ u , v \}$ (and $u \wedge v:= \min \{ u,v \} $). 
For $N \in \Nat$ and a vector $w\in \R^N$ and a real-valued function $f$ we define 
the vector $f(w)$ as $f(w)  := (f(w_1)  , \dots , f( w_N))^T $.


    \section{The noisy case}\label{noisy.section}
     
  This section studies the noisy model
  $$ Y = X \beta^0  +\epsilon , $$
  where $Y$ is an $n$-vector of observations and 
 with $\epsilon= (\epsilon_1 , \ldots , \epsilon_n)^T$ containing  i.i.d. ${\cal N} (0,1/n ) $-distributed noise  variables.
  We will compare the noisy Lasso
 $$ \hat \beta := \arg \min_{\beta \in \R^p } \biggl \{
 \| Y - X\beta \|_2^2 + 2 \lambda \| \beta \|_1 \biggr \} $$
 with the noiseless Lasso
 $$ \beta^* := \arg \min_{\beta \in \R^p}  \biggl \{  \| X (\beta - \beta^0) \|_2^2 + 
 2 \lambda \| \beta \|_1\biggr \} . $$
 We show in the next two theorems that under certain conditions on the design the
 ``bias" $ \| X (\beta^* - \beta^0) \|_2$ is of larger order (in probability) than the ``estimation error"
 $ \| X (\hat \beta - \beta^* ) \|_2$ (where ``bias" and ``estimation error" are here to be understood
 in generic terms). By the triangle inequality
 $$ \| X ( \hat \beta - \beta^0 )\|_2 \ge  \underbrace{ \| X ( \beta^* - \beta^0) \|_2}_{\rm ``bias"}  -
\underbrace{ \| X ( \hat \beta - \beta^* ) \|_2 }_{\rm ``estimation \ error"}  $$
this implies a high probability lower bound for  the prediction error $\| X ( \hat \beta - \beta^0) \|_2$
of the noisy Lasso in terms of the prediction error $\| X ( \beta^* - \beta^0) \|_2$ of the noiseless Lasso. 
 
  \begin{theorem}\label{silly.theorem} Let $\| X_j \|_2\le 1$ for all $j$, and
 let $0<\alpha < 1$ and $0 < \alpha_1 < 1$ be fixed and 
 $\lambda_0 := \sqrt {2 \log (2p/\alpha) /n} $. Let $\eta \lambda > \lambda_0$ for some
 $0 \le  \eta < 1$.  Then with probability at least
 $1- \alpha - \alpha_1$
 $$ \| X ({\hat \beta} -{\beta^*} ) \|_2 \le \sqrt {  \Lambda_{\rm max} 
 (\hat \Sigma) \over n \lambda^2(1- \eta)^2  } { \|  X( \beta^* - \beta^0) \|_2  }  +
 \sqrt {2 \log (1/\alpha_1)\over n  } .$$
  \end{theorem}
  
  {\bf Asymptotics} We see we may choose $\lambda \asymp \sqrt {\log p / n} $.
  Then, for $p \rightarrow \infty$ and $\Lambda_{\rm max} ( \hat \Sigma ) = 
  {\mathcal O} (1)$  we get
  $$ \| X ({\hat \beta} -{\beta^*} ) \|_2  = o_{\PP} (1)  \| X ({\beta^*} - \beta^0) \|_2 + {\mathcal O}_{\PP} (1/ \sqrt n) . $$
  
  In general the largest eigenvalue $\Lambda_{\rm max} ( \hat \Sigma) $ may be large,
  and may be hard to control, for example when the Gram matrix $\hat \Sigma$ comes from
  random design. 
  We now let $\Sigma_0$ be some approximation of $\hat \Sigma$, for example
  a population version $\EE \hat \Sigma$ of $\Sigma_0$ in the case of random design.
  
  We use the notation $\| \hat \Sigma - \Sigma_0 \|_{\infty} := \max_{j,k} | \hat \Sigma_{j,k} - \Sigma_{0,j,k} | $. 
  
   \begin{theorem}\label{sillybis.theorem} Let $\| X_j \|_2\le 1$ for all $j$, and let $0<\alpha < 1$ and $0 < \alpha_1 < 1$ be fixed and 
 $\lambda_0 := \sqrt {2 \log (2p/\alpha) /n} $. Let $\eta \lambda > \lambda_0$
 for some $0 \le  \eta < 1$.  
 Suppose that
 \begin{eqnarray}\label{small.equation}
 \xi:= \| \hat \Sigma - \Sigma_0 \|_{\infty} \| \beta^* - \beta^0 \|_1 < \lambda (1- \eta) .
 \end{eqnarray}
 Then with probability at least
 $1- \alpha - \alpha_1$
 $$  \| X ({\hat \beta} -{\beta^*} ) \|_2 \le { \Lambda_{\rm max}^{1/2}  ( \Sigma_0 )  \biggl ( \| X( \beta^* - \beta^0)  \|_2^2  + \xi  \| \beta^* - \beta^0 \|_1 \biggr )^{1/2}  \over 
  \biggl (  \lambda (1- \eta )- \xi  \biggr )} 
  +
 \sqrt {2 \log (1/\alpha_1) \over n  }.
 $$ 
  \end{theorem}
  
  Condition (\ref{small.equation}) is a condition requiring the $\ell_1$-error of $\beta^*$ to be small. In an asymptotic
  setup, it typically needs sparsity $s_0$ of small order $\sqrt {n / \log p }$. However,  
  in the case of Gaussian random design for example and $\Sigma_0= \EE \hat \Sigma$ one may apply more
  careful bounds to prove a result that does not require such sparsity conditions.

\section{Upper bounds}\label{upper-bound.section}
There are several upper bounds in the literature. The one we will mainly apply
is along the lines of Theorem 6.1 in \cite{BvdG2011}, with some refinements.
The result is given in Lemma \ref{upperbound1.lemma}.
There are however more general bounds in literature, in particular {\it sharp oracle bounds}
as in \cite{koltchinskii2011nuclear} (see also 
\cite{giraud2014introduction}, Theorem 4.1 or \cite{vdG2016}, Theorem 2.2).
We present these in Lemma \ref{upperbound.lemma}.

The upper bounds follow from the KKT (Karush-Kuhn-Tucker) conditions
$$\hat \Sigma  (\beta^* - \beta^0) + \lambda z^* =0 .$$
Here $z^* \in \partial \| \beta^* \|_1$ with $\partial \| \beta\|_1 $ the sub-differential
of the mapping $\beta \mapsto \| \beta \|_1$, $\beta \in \R^p$. In other words
$\beta^{*T}z^*= \| \beta^* \|_1$ and $\| z^*\|_{\infty} \le 1 $. 

Here are the  upper bounds for the prediction error we will use. They 
include upper bounds for $\| \beta^* \|_1$ and $\| \beta_{-S_0}^* \|_1 $.

\begin{lemma}\label{upperbound1.lemma} It holds that
$$  \| X ( \beta^* - \beta^0) \|_2^2 + \lambda \| \beta^*\|_1 \le \lambda \| \beta^0 \|_1 , $$
and
$$ \| X ( \beta^* - \beta^0) \|_2^2 + 2 \lambda \| \beta_{-S_0}^*\|_1\le
{\lambda^2 s_0 \over \hat \phi^2 (S_0) } . $$
\end{lemma}

The next lemma contains the more general sharp oracle inequalities for the
prediction error. 

\begin{lemma} \label{upperbound.lemma} 
The prediction error  $\| X (\beta^* - \beta^0 ) \|_2^2 $ satisfies the bound
 \begin{eqnarray*}
  \| X (\beta^* - \beta^0 ) \|_2^2  \le  \hat {\cal U} ( \beta^0) ,
  \end{eqnarray*}
  where
  $$ 
    \hat {\cal U}( \beta^0) = \min \{\hat {\cal U}_I( \beta^0), \hat {\cal U}_{II}( \beta^0), \hat {\cal U}_{III}( \beta^0)\} 
   $$
   with 
  \begin{eqnarray*}
 & \ &  \hat {\cal U}_{I}( \beta^0):= { \lambda^2 s_0 \over \hat \phi^2 (S_0) } \wedge  \lambda \| \beta^0 \|_1 \\
 & \ & 
 \hat {\cal U}_{II}( \beta^0)\\
 & \ &  \ \ := \min_S \left \{ \biggr (\sqrt {\lambda^2 |S| \over 4 \hat \phi^2 (S)} + \sqrt {{ \lambda^2 |S| \over   4 \hat \phi^2 ( S) }
  + \lambda \| \beta_{-S}^0 \|_1  } 
   \biggr  ) ^2  \vee 2 \lambda  \| \beta_{-S}^0 \|_1 \right \} \\
  & \ & \hat {\cal U}_{III}( \beta^0)  \\ & \ & \ \  := 
  \min_{S} \min_{\beta}  \left \{ \biggl (
 \| X(\beta - \beta^0) \|_2^2 + { \lambda^2 |S| \over  \hat \phi^2 (S) } + 2 \lambda \| \beta_{-S} \|_1 
 \biggr ) \vee 4 \lambda \| \beta_{-S} \|_1 \right \} 
  \end{eqnarray*}    . 
 \end{lemma}
 
 Clearly, if \\
 $\circ$ the minimum over $S$ in the definition of $\hat {\cal U}_{II} (\beta^0)$ is attained
 in $S_0$,\\
 $\circ$ the minimum over $(S,\beta)$ in the definition of
 $\hat {\cal U}_{III} (\beta^0)$ is attained in $(S_0, \beta^0)$,\\
 and\\
 $\circ$ $\lambda^2 s_0 / \hat \phi^2 (S_0) \le \lambda \| \beta^0 \|_1 $,
 then
 $$\hat {\cal U}_{I} (\beta^0) = \hat {\cal U}_{II} (\beta^0)= \hat {\cal U}_{III} (\beta^0)=
{  \lambda^2  s_0 \over \hat \phi^2 (S_0) } . $$
This will be the case in most of the examples we consider in this paper, that is, we do not explore
the power of the sharp oracle inequalities of Lemma \ref{upperbound.lemma}. Instead,
we mainly compare exact results for the (penalized) prediction error with the bounds of Lemma \ref{upperbound1.lemma}.

 \begin{remark} \label{U.remark}
 Clearly, Lemma \ref{upperbound1.lemma} implies the bound $\hat {\cal U}_{I}( \beta^0)$. 
 Further, by restricting $S$  in the minimization giving $\hat {\cal U}_{II}( \beta^0)$ to 
$S \in  \{ S_0 ,  \emptyset \} $ one sees
 $$ \hat {\cal U}_{II}( \beta^0) \le { \lambda^2 s_0 \over \hat \phi^2 (S_0) } \vee 2 \lambda  \| \beta^0 \|_1 . $$
 In other words,  up to a factor ``2",  the bound $\hat {\cal U}_{II}( \beta^0)$ improves upon 
 $ \hat {\cal U}_{I}( \beta^0)$. Similarly, taking  $\beta = \beta^0$ in the 
 minimization giving $\hat {\cal U}_{III}( \beta^0)$ one finds
 $$  \hat {\cal U}_{III}( \beta^0)\le \biggl (  { \lambda^2 |S| \over  \hat \phi^2 (S) } + 2 \lambda \| \beta_{-S}^0 \|_1 
 \biggr ) \vee 4 \lambda \| \beta_{-S}^0 \|_1  ,$$
 that is, up to a factor ``2",  $   \hat {\cal U}_{III}( \beta^0)$ improves upon  $\hat {\cal U}_{II}( \beta^0)$. 
 Note also that 
 \begin{eqnarray} \label{project.equation}
  \hat {\cal U}_{III}( \beta^0)\le 
\min_{S}  \left \{ 
 \| X(b_S - \beta^0) \|_2^2 + { \lambda^2 |S| \over  \hat \phi^2 (S) } 
 \right \} 
 \end{eqnarray}
 where (for every set $S $) $X b_S $ is the projection of $X \beta^0$ on the space spanned by
 $ \{ X_j \}_{j \in S} $.
 \end{remark}

\section{Some considerations about the design}\label{design.section}

 \begin{definition} We say that $X$ has normalized columns if for any $j$ it holds that
  $ \| X_j \|_2   = 1 $. We then  call the design normalized. 
   \end{definition}
   
    \begin{definition} We say that $X$ has no aligned columns if for any $j \not= k$, and any constant $b $ it holds that
 $X_j \not= b X_k $.
 \end{definition}
 
 \begin{definition} We say that $X$ is a fair design if it is normalized and has no aligned columns.
  \end{definition}
  
  The reason for requiring normalized design is that when the columns in $X$ have different lengths,
  say the length of the first column $X_1$ is much smaller than that of the others, then in effect
  the first variable gets a heavy penalty as compared to the others. By taking $\| X_1 \|_2$ extremely 
  small, one can force the Lasso to choose $\beta_1^*$ extremely small, thus creating an unfair
  situation. 
  
  With normalized design, no aligned columns means that  $X_j \not= \pm X_k $ for all $j \not= k$.
  
  As we will see, one of the reasons why in the rest of the paper we assume that there are at least
  two active variables is the following:
  
  \begin{lemma} \label{1fair.lemma} There is no fair design such that $\hat \phi(\{1\})=0$.
   \end{lemma}

\section{Assumption about the first two variables}\label{first2.section}
In what follows we consider throughout the case where $\beta_1^0 \ge \beta_2^0 > 0$ so that the
first two variables are among the active ones. Moreover, we assume
$$(X_1 , X_2)^T (X_1 , X_2)    = \begin{pmatrix} 1 & - \hat \rho \cr - \hat \rho & 1 \cr  \end{pmatrix}, $$
where $0 < \hat \rho = -X^T X_2   < 1 $ is minus the inner product between $X_1$ and $X_2$.
Although we do not insist that $X_1$ and/or $X_2$ are centered, we sometimes refer
to $-\hat \rho$ as the correlation between $X_1$ and $X_2$. 
The negative correlation is to be seen 
in relation with both $\beta_1^0$ and
$ \beta_2^0$ positive. It is 
so to speak the more difficult case for the Lasso. 

Throughout the paper, we set
$$\hat \varphi^2 := 1 - \hat \rho. $$
Fair design as defined in the previous section is related to using the penalty $\lambda \| \beta |_1$ with equal
weights for all coefficients. But  linear combinations of the columns in $X$ are of
course generally not normalized. We obviously have for example
$ \| X_1 + X_2 \|_2^2 = 2\hat \varphi^2 $ which is less than 1 when $\hat \phi^2 < 1/2$. 
As we will see this is roughly the main ingredient when constructing exact results depending
on compatibility constants. 

\section{Results for $p=s_0=2$} \label{p=s0=2.section}
In this section $p$ equals $2$ so that $X= (X_1 , X_2) $. One may argue that this is not
exactly a high-dimensional situation (for which the Lasso is designed)
and therefore of limited interest.  However, lower bounds for
the low-dimensional situation can easily be extended to higher dimensions (trivially for example,
by adding inactive variables orthogonal to the active ones, see Section \ref{trivial.section}).
If the irrepresentable condition holds, the Lasso
will not select inactive variables (see  \cite{ZY07}) which brings
us back in a lower-dimensional situation. 
Lemmas \ref{phi=rho.lemma} and 
 \ref{2comp>0.lemma}
are examples where the Lasso ignores inactive variables
that are correlated with the active ones.

\begin{lemma} \label{hatphi1.lemma} We have
$$ \hat \phi^2 (\{1\}) = 1- \hat \rho^2 =   \hat \varphi^2  (2-\hat \varphi^2)  . $$
Moreover
$$ \hat \phi^2 ( S_0 ) \ \  = \hat \varphi^2   , \ \hat \Gamma^2 (S_0) = {2 \over \hat \varphi^2 } . $$
\end{lemma}

In the case considered here ($p=2$) the minimal eigenvalue $\Lambda_{\rm min } ( \hat \Sigma ) $  of the
Gram matrix $\hat \Sigma$ is
$$ \Lambda_{\rm min} (\hat \Sigma) = 1- \hat \rho =\hat \varphi^2  .$$
Thus, the compatibility constant $\hat \phi^2 (S_0)$ is just another expression for 
this minimal eigenvalue.
Lemma \ref{eigenvalue.lemma} gives an example in  a higher-dimensional case, where the compatibility constant can be
(much) larger than  $\Lambda_{\rm min} (\hat \Sigma)$, and in fact also (much) larger than the
restricted eigenvalue as defined in \cite{bickel2009simultaneous}.

\begin{lemma}\label{p=2-exact.lemma} 
Consider the following three cases:
\begin{eqnarray*}
{\rm Case \ 1: } & \ & { \lambda / \hat \varphi^2  } \le  \beta_2^0 \\
{\rm Case \ 2:} & \ &  \beta_2^0 \le { \lambda / \hat \varphi^2 } \le  \beta_2^0 +
{ (\beta_1^0 - \beta_2^0) /  \hat \varphi^2  }\\
{\rm Case \ 3:} &\ & { \lambda /\hat \varphi^2  } \ge   \beta_2^0 +
{ (\beta_1^0 - \beta_2^0 ) /  \hat \varphi^2  }  . 
\end{eqnarray*} 
Then we have
$$ \| X (\beta^* - \beta^0) \|_2^2 = \begin{cases} {2 \lambda^2 /  \hat \varphi^2  } & \ {\rm in \ Case \ 1} \cr
\hat \varphi^2 (2- \hat \varphi^2) (\beta_2^0)^2+ \lambda^2 
  &  \ {\rm in \ Case \ 2} \cr 
 \| X \beta^0 \|_2^2  &  \ {\rm in \ Case \ 3}
 \end{cases} \ \  $$
 and
  $$ \beta^*  = \begin{cases}   \begin{pmatrix} \beta_1^0 - \lambda/ \hat \varphi^2  \cr
  \beta_2^0 - \lambda/ \hat \varphi^2 \cr 
  \end{pmatrix} & \ {\rm in \ Case \ 1} \cr
  \begin{pmatrix} \beta_1^0-
 (1- \hat \varphi^2 ) \beta_2^0 -\lambda \cr  0\cr  \end{pmatrix}
 &  \ {\rm in \ Case \ 2} \cr  \ \ \begin{pmatrix} 0\cr 0 \cr \end{pmatrix}
 &  \ {\rm in \ Case \ 3}
 \end{cases} \ \ . $$
\end{lemma}

\begin{corollary} Lemma \ref{p=2-exact.lemma} reveals that in ${\rm Case \ 1 }$ 
$$ \| X ( \beta^* - \beta^0 )\|_2^2 + \lambda \| \beta^* \|_1 =
\lambda \| \beta^0 \|_1 , $$
and, invoking Lemma \ref{hatphi1.lemma},
  \begin{eqnarray*}
 \| X (\beta^* - \beta^0)  \|_2^2   =  { \lambda^2 s_0  \over \hat \phi^2 (S_0)  }    . 
 \end{eqnarray*}
  This corresponds exactly to the bounds in Lemma \ref{upperbound1.lemma}. 
\end{corollary} 

\begin{corollary}\label{p=2-exact.corollary}
It may be of interest to consider the intersection of the  cases in Lemma \ref{p=2-exact.lemma}. We see
that
$$ \| X (\beta^* - \beta^0) \|_2^2  $$ $$ = \begin{cases} {2 \lambda \beta_2^0  } & {\rm in \ Case \ 1 \cap 2 }:
 \  \beta_2^0 ={ \lambda /  \hat \varphi^2 } 
\cr
2 \lambda \beta_1^0 - ((\beta_1^{0})^2 - (\beta_2^{0})^2  ) & {\rm in\ Case \ 2 \cap 3}:\ 
 { \lambda / \hat \varphi^2  } =  \beta_2^0 +
{ (\beta_1^0 - \beta_2^0) / \hat \varphi^2  }  \cr 
 \end{cases} . $$
 Thus, the bound $\hat {\cal U}_{II} (\beta^0)$ in Lemma
 \ref{upperbound.lemma} is tight in ${\rm Case \ 1 \cap 2 }$. 
\end{corollary}

\begin{corollary}\label{equal-beta.corollary}
When $\beta_1^0 = \beta_2^0 $, 
the union of cases gives
$$  \| X ( \beta^* - \beta^0) \|_2^2   = \begin{cases} 2 \lambda^2 / \hat \varphi^2   & {\rm \ in \  Case \ 1 \cup 2 }:\ 
\lambda / \hat \varphi^2  \le  \beta_2^0 \ \cr 2 \hat \varphi^2  (\beta_2^0)^2 & {\rm \ in \ Case \ 2 \cup 3 }:\ 
\lambda / \hat \varphi^2  \ge \beta_2^0\cr \end{cases}. $$
\end{corollary}

 \begin{remark} On may verify that ${\rm Case \ 2} $ has
$$ \| X ( \beta^* - \beta^0) \|_2^2 =  \| X( b_{\{1 \}} - \beta^0 ) \|_2^2 + \lambda^2 , $$
where $Xb_{\{ 1 \}} $ is the projection of $X \beta^0 $ on $X_1$. 
This can be compared with (\ref{project.equation}) (following from
 ${\cal U}_{III} (\beta^0)$ defined
in Lemma \ref{upperbound.lemma}) in Remark \ref{U.remark}. 
  \end{remark}
  
  \begin{remark} \label{Case-3.remark} In ${\rm Case \ 3}$, we have
$ {( \lambda - (\beta_1^0  - \beta_2^0) )/ \hat \varphi^2  } \ge  \beta_2^0 >0 $. This
implies $\beta_1^0 - \beta_2^0 <  \lambda$. Note moreover that this case illustrates that the bound
(\ref{project.equation}) in Remark \ref{U.remark} (and hence ${\cal U}_{III} (\beta^0)$ defined
in Lemma \ref{upperbound.lemma}) can be tight.
\end{remark} 

\begin{remark} \label{extreme-case.remark} The case $\hat \rho =  1$ is not treated in Lemma
\ref{p=2-exact.lemma}.  It corresponds to Case 2 with $\hat \varphi^2  \downarrow 0$.
\end{remark}

 \section{Results for $p=s_0 = 2N$}\label{p=s0=2N.section}
 The results of the previous section are easily extended to a larger active set $S_0$. We assume
 $S_0 = \{ 1 , 2 , \ldots , s_0 \}$ with $s_0$ even, say $s_0 = 2N$ (with $N \in \Nat $ and $2N \le n$). Moreover we again assume
 $p=s_0$. Then
 $$ X= (X_1 , X_2 , \ldots , X_{2N-1}, X_{2N}  ) .$$
  We split the design into $N$ matrices of dimension $n \times 2 $.
 
 \begin{lemma}\label{eigenvalue.lemma} Consider fair design with (for $k\in \{ 1, \ldots , N\}$)
$(X_{2k-1}, X_{2k})$ orthogonal to the space spanned by 
 the remaining columns.  Assume that  $\hat \rho_k := -X_{2k-1}^T X_{2k}>0$  
 and  write $\hat \varphi_k^2  := 1- \hat \rho_k $ for all $k\in \{ 1, \ldots , N \} $.
Then 
$$\Lambda_{\rm min} (\hat \Sigma) =  \min_{k} \hat \varphi_k^2  $$
and
$$ \hat \phi^2 ( S_0) = { N \over \| 1/ \hat \varphi^2 \|_1 } \ge  \Lambda_{\rm min} ( \hat \Sigma)  ,  \ \hat \Gamma^2 (S_0) =2  \| 1 / \hat \varphi^2 \|_1 . $$
Moreover, for ${\cal S}= \{ 2, 4, \ldots , 2N \}$
$$ \hat \phi^2 ({\cal S})=  {  N \over \| (1- \hat \rho^2)^{-1}  \|_1}    .$$
 \end{lemma}

\begin{remark}The restricted eigenvalue (\cite{bickel2009simultaneous}) is defined as
$$ \hat \kappa^2 (S) = \min \biggl \{ { \| X \beta_S - X \beta_{-S} \|_2^2  \over
\| \beta_S \|_2^2 } :\ \| \beta_{-S} \|_1 \le \| \beta_S \|_1 \biggr \} . $$
In the case we are considering in this section, where $S_0= \{1 , \ldots , p \} $,
one obviously has $ \hat \kappa^2 (S_0)= \Lambda_{\rm min} ( \hat \Sigma)$.
Therefore,  in the situation of Lemma \ref{eigenvalue.lemma}
$ \hat \kappa^2 (S_0) \le \hat \phi^2 (S_0)$ and the difference can be substantial.

\end{remark}

 The next lemma is again an illustration of the tightness of the upper bounds
 in Lemma \ref{upperbound1.lemma}. 
 
 \begin{lemma} \label{UII-IIIbis.lemma} Consider design as in Lemma \ref{eigenvalue.lemma}.
 Suppose that for all
 $k$, $\beta_{2k-1}^0 \ge \beta_{2k}^0 \ge \lambda /\hat \varphi_k^2  $. Then
 $$  \| X (\beta^* - \beta^0) \|_2^2 ={\lambda^2 s_0\over \hat \phi^2 (S_0) } = 
 { 2 \lambda^2  \| 1/ \hat \varphi^2  \|_1 } $$
 and
 $$ \| X (\beta^* - \beta^0) \|_2^2 + \lambda \| \beta^* \|_1 = \lambda \| \beta^0 \|_1  . $$
     \end{lemma}

   \begin{remark}
     For the special case of Lemma \ref{UII-IIIbis.lemma} with equality
     $ \beta_{2k}^0 =\lambda/ \hat \varphi_k^2   $ for all $k$, have 
    for ${\cal S} := \{2 , 4, \ldots , 2N \} $
    $$  \| X (\beta^* - \beta^0) \|_2^2
   = 
 2 \lambda \| \beta_{-{\cal S} } \|_1 ,$$
 showing tightness of $\hat {\cal U}_{II} (\beta^0) $.
 \end{remark}

   \section{A trivial extension to $m_0 >0$}\label{trivial.section}
   Recall that $X_{-S_0}$ contains the $m_0 := |S_0^c|$ inactive variables. If these are orthogonal
   to the active ones the results are trivially as for the case $m_0=0$.
   As an example, let us take $s_0=2$.
   
     \begin{lemma} \label{trivial.lemma} Let again $S_0 = \{1 , 2\}$  and  suppose that $X_{S_0}^T X_{-S_0} = 0 $. Then
  $$ \hat \phi^2 (S_0) = \hat \varphi^2  $$
  and for $\beta_1^0 \ge \beta_2^0 \ge \lambda / \hat \varphi^2 $, 
  $$ \| X (\beta^*- \beta^0 ) \|_2^2 + \lambda \| \beta^* \|_1 = \lambda \| \beta^0 \|_1  $$
  and
  $$ \| X (\beta^*- \beta^0 ) \|_2^2 = { \lambda^2 s_0  \over \hat \phi^2 (S_0) }=
   {2 \lambda^2 \over \hat \varphi^2 } . $$
   \end{lemma}

   By the same argument, one may always extend in what follows the non-active set
   with variables that are orthogonal the ones considered.

\section{A result for $s_0 = 2$, $m_0=1$}\label{s0=2-m0=1.section}
We now add one inactive variable, that is we take
$S_0= \{1 , 2 \} $ and $S_0^c= \{ 3 \} $. 

\begin{lemma} \label{s0=2-m0=1.lemma} Suppose that
$$ X_3 = C (X_1 + X_2)/2 + U $$
where $C$ is a constant satisfying $C > 1 $ and $C^2 \hat \varphi^2 /2< 1 $,
and $U$  is a vector with $U^T (X_1 , X_2) =0$.
Define
$$ \hat \tau^2 := 1-  C^2 \hat \varphi^2/2 . $$
Then
$$ \hat \phi^2 (S_0) = \hat \varphi^2 \hat \tau^2 , \ \hat \Gamma^2 (S_0):= { s_0 \over \hat \phi^2 (S_0) } =
{ 2 \over \hat \varphi^2 } + {  C^2  \over \hat \tau^2 } .$$
For $\beta_1^0 \ge \beta_2^0 \ge \lambda /\hat \varphi^2 + \lambda C (C-1) / (2\hat \tau^2)$
we have
$$ \| X(\beta^* -\beta^0 )\|_2^2 + 2\lambda \| \beta_{-S_0}^* \|_1 =
{  \lambda^2 s_0 \over \hat \phi^2 (S_0)} - {\lambda^2 \over \hat \tau^2} . $$
\end{lemma}

The above lemma shows that the second upper bound of Lemma \ref{upperbound1.lemma}
is   a term $\lambda^2 / \hat \tau^2 $ too large.
However, this term can be small. An example is given in the next corollary.

\begin{corollary} Take in Lemma \ref{s0=2-m0=1.lemma} the constant $C=2$.
Then for $\hat \varphi^2 < 1/2$
$$ \hat \phi^2 (S_0)= \hat \varphi^2 (1- 2\hat \varphi^2) , $$
and so
$$ \hat \Gamma^2 (S_0) = {  2  \over \hat \varphi^2} + { 4   \over 1- 2 \hat \varphi^2 } , $$
and for $\beta_1^0 > \beta_2^0 \ge \lambda (1- \hat \varphi^2 )/ (\hat \varphi^2 (1- 2 \hat \varphi^2 ) $
we have
$$ \| X(\beta^* -\beta^0) \|_2^2 + 2\lambda \| \beta_{-S_0}^* \|_1 = 
{ 2 \lambda^2 \over \hat \varphi^2} + {3 \lambda^2 \over 1- 2 \hat \varphi^2 } . $$
In other words, the bound in lemma \ref{upperbound1.lemma} has a factor
``4" whereas the exact result has a factor ``3". For $\hat \varphi^2 \downarrow 0$ we see
that the upper bound is asymptotically tight, as then
${ 2 \lambda^2 / \hat \varphi^2}$ is the leading term. Conversely,
for $\hat \varphi^2 \uparrow 1/2$ the upper bound is asymptotically a factor $4/3$ too large. 
\end{corollary}

\section{A result for $s_0 = 2N$, $m_0=1$}\label{s0=2N-m0=1.section}
We have seen in the previous section that the upper bound of Lemma \ref{upperbound1.lemma}
can be off, for example by a factor $4/3$ asymptotically. The question arises whether in a generalized
setting
this factor increases when $s_0$ increases. If this is not the case, the non-tightness
of the bound is really only a matter of constants. In this section
we show in an example that the gap between the upper bounds of Lemma
\ref{upperbound1.lemma} and the exact bound does not depend on $s_0$.

\begin{lemma}\label{s0=2N-m0=1.lemma} 
Let $S_0 = \{1 , \ldots , 2N\} $, $S_0^c= \{ 2N+1 \} $ and
$$ (X_{2k-1}, X_{2k} )^T(X_{2k-1}, X_{2k} )  := 
\begin{pmatrix} 1 & - \hat \rho_k \cr - \hat \rho_k  & 1 \cr \end{pmatrix}, \ k=1 , \ldots , N , $$
where each $\hat \rho_k$ is between 0 and 1. Then we define
$\hat \varphi_k^2 :=1 - \hat \rho_k $, $k=1 , \ldots , N$. 
Further, assume that $(X_{2k-1} , X_{2k } ) $ is orthogonal
to $\{ X_j \}_{j \in S_0\backslash  \{ 2k-1 , 2k  \},}  $ for all $k$. Let
$$X_{2N+1} = C \sum_{j=1}^{2N} X_j / s_0 +U $$ where $ C >  1$, 
$C^2 \sum_{k=1}^N 2 \hat \varphi_k^2 / s_0^2 < 1 $ and $U$ is orthogonal
to $\{ X_j \}_{j \in S_0}$. 
Write $\hat \tau^2 :=1- C^2 \sum_{k=1}^N 2 \hat \varphi_k^2 / s_0^2 $.
Then for $\hat \varphi^2 = ( \hat \varphi_1^2 , \cdots , \hat \varphi_k^2)$
$$ \hat \Gamma^2 (S_0) = 2 \| 1/ \hat \varphi^2 \|_1 + {C^2 \over \hat \tau^2 } . $$
Moreover, for $\beta_{2k-1}^0 \ge \beta_{2k}^0 \ge  \lambda / \hat \varphi_k^2 + \lambda C(C -1) / (s_0 \hat \tau^2) $, 
$k=1 , \ldots , N$, we have
$$ \| X (\beta^* - \beta^0 )\|_2^2 + 2 \lambda \| \beta_{-S_0}^* \|_1 =
\lambda^2 \hat \Gamma^2 (S_0)  - {\lambda^2 / \hat \tau^2 } . $$
\end{lemma}

\begin{corollary}\label{s0=2N-m0=1.corollary} When $\hat \varphi_1^2 = \cdots = \hat \varphi_N ^2:= \hat \varphi_0^2$
 (say)
 in Lemma \ref{s0=2N-m0=1.lemma} and $C= 2 $ one gets 
 $$ \hat \Gamma^2 (S_0)= { s_0 \over \hat \varphi_0^2} + { 4 \over \hat \tau^2 } , $$
 with $\hat \tau^2 = 1- 4 \hat \varphi_0^2 / s_0$.  For
 $\beta_{2k-1}^0 \ge \beta_{2k}^0 \ge \lambda / \hat \varphi_k^2+ 2\lambda / (\hat \tau^2 s_0) $
 for all $k$, 
 we get 
 $$ \| X (\beta^* - \beta^0 )\|_2^2 + 2 \lambda \| \beta_{-S_0}^* \|_1 =
{ \lambda^2 s_0  \over \hat \varphi_0^2 }  - {\lambda^2  \over \hat \tau^2 } . $$
So with $\hat \varphi_0^2$ kept fixed  the gap of Lemma \ref{upperbound1.lemma} decreases with $s_0$. 
 \end{corollary}

\section{ A result for $s_0 = 2 $ and $m_0$ possibly large}
\label{s0=2-m0.section}
We now set $S_0 = \{1 , 2 \}$ and
$S_0^c:= \{ 3 , \ldots, 2+ m_0 \}$ where $m_0$ is possibly large
(in an asymptotic sense it may be of order $1/\lambda$ say). 

\begin{lemma} \label{s0=2-m0.lemma} Suppose
$$ X_{2+k} = C_k ( X_1 + X_2)/2 + U_k ,\  k=1, \ldots , m_0 , $$
where, for $k=1 , \ldots , m_0$, the constant $C_k $ has $C_k >1 $ but
$ C_k^2 \hat \varphi^2 /2 < 1 $, and where the vector  $U_k$ is orthogonal to 
$\{ X_1 , X_2 , 
\{ U_j \}_{j \not= k } \} $. \\
Let for each $k  \in \{ 1 , \ldots , m_0 \} $, the constant
$\hat \tau_k^2$ be given by $\hat \tau_k^2 = 1-  C_k^2 \hat \varphi^2/2 $.
Then
$$ \hat \Gamma^2 (S_0)= { 2 \over \hat \varphi^2} + \sum_{k=1}^{m_0} 
{  C_k^2 \over \hat \tau_k^2 } . $$
Moreover, if $\beta_1^0 \ge \beta_2^0 \ge \lambda / \hat \varphi^2 + \lambda \sum_{k=1}^{m_0}
C_k ( C_k-1) / (2\hat \tau_k^2 )$, it holds that
$$ \| X ( \beta^* - \beta^0 ) \|_2^2 +
2 \lambda \| \beta_{-S_0}^* \|_1 = \lambda^2 \hat \Gamma (S_0) - \lambda^2 \|  1/ \hat \tau^2 \|_1 $$
\end{lemma}

\begin{corollary}\label{s0=2-m0.corollary} If we take $C_k= 2$ for all $k \in \{1 , \ldots , m_0 \}$
we obtain
$$ \hat \Gamma^2 (S_0) = {2 \over \hat \varphi^2 } + { 4 m_0 \over 1- 2 \hat \varphi^2 } , $$
and
$$ \| X ( \beta^* - \beta^0) \|_2^2 + 2 \lambda \| \beta_{-S_0} \|_1 =
{ 2 \lambda^2 \over \hat \varphi^2 } + { 3 \lambda^2 m_0 \over 1- 2 \hat \varphi^2 } . $$
The upper bound of Lemma \ref{upperbound1.lemma}
is off no more than a factor $4/3$. 

\end{corollary}

  \section{Some results for $s_0=m_0 =2$}
  \label{s0=2=m0.section}
  
  In this section, the active set is again $S_0 = \{1 , 2 \}$ and the non-active one is $S_0^c = \{ 3, 4 \}$. Thus, both $ s_0 $ and $m_0:= p-s_0$ are equal to 2. 
  
  In Section \ref{s0=2-m0=1.section}, we have seen that the upper bound of Lemma \ref{upperbound1.lemma}
  can be  too large, but that the gap is small 
  when the main term is due to highly negatively correlated active variables.
  In this section,  we consider  first a
  setup similar to the one in Section \ref{s0=2-m0=1.section}. Again, the upper bounds
  are not tight but the gap can be small. Unlike the previous section,
  the main terms in the  bound in this section are now not necessarily determined by
  the negative correlations in the active set.

    \begin{lemma}\label{goodcomp.lemma} Let
    $$ X_3 = C(X_1 + X_2 )/2+ U+V , \ X_4 = C(X_1 + X_2)/2 + U-V , $$
    where $C>1$, $ C^2 \hat \varphi^2 /2< 1 $,  $U^T X_{S_0} = V^T X_{S_0} =0 $ and $U^T V =0$. Set
    $$ \hat \tau^2 :=  U^T U  $$
    where $0< \hat \tau^2 <  1 -  C^2\hat \varphi^2/2 $. 
    Then
    $$ \hat \phi^2 (S_0)= {\hat \varphi^2 \hat \tau^2 \over C^2 \hat \varphi^2/2  + \hat \tau^2 } , \
    \hat \Gamma^2 (S_0) = { 2 \over \hat \varphi^2 } + { C^2 \over \hat \tau^2} .   $$
    Let $\beta_1^0 \ge \beta_2^0 \ge \lambda/ \hat \varphi^2 + \lambda C(C-1) / (2\hat \tau^2 )$.
 Then
 $$ \| X ( \beta^* - \beta^0) \|_2^2 + 2 \lambda \| \beta_{-S_0^*} \|_1 = \lambda^2 \hat \Gamma^2 (S_0)  - { \lambda^2 \over \hat \tau^2}   . $$
  \end{lemma}

    We can also have a look what happens if in the above lemma, we let
   $\hat \tau^2  =0$
    instead of $>0$. Then the compatibility constant  $\hat \phi^2 (S_0) $ is zero.
    In this case, the prediction error 
    $\| X ( \beta^* - \beta^0 ) \|_2^2$ is in a sense still under control, but the penalized prediction
    error  $\| X ( \beta^* - \beta^0 ) \|_2^2 + 2 \lambda \| \beta_{-S_0} \|_1 $ can show the ``slow rate".
    
    \begin{lemma}\label{goodlasso2.lemma}
    Let
    $$ X_3 = C(X_1 + X_2 )/2 +V , \ X_4 = C(X_1 + X_2)/2 -V , $$
    where $C>1$, $ C^2 \hat \varphi^2/2 < 1 $ and $ V^T X_{S_0} =0 $.
        Then
    $$ \hat \phi^2 (S_0)=0 . $$
    Moreover when $\beta_1^0 \ge  \beta_2^0 \ge \lambda / \hat \varphi^2 $ we find
    $$ \| X ( \beta^* - \beta^0) \|_2^2 = {2 \lambda^2 \over \hat \varphi^2 }, $$
    $$ \| X ( \beta^* - \beta^0) \|_2^2 + 2 \lambda \| \beta_{-S_0}^* \|_1 = 
    {4 \lambda \beta_2^0  \over C}  - { 2 \lambda^2 \over \hat\varphi^2 } \biggr ( {2 \over C}-1 \biggl ) .  $$
     \end{lemma}
     
     Note that if in the above lemma $C=2$ we arrive at the bound
      $$ \| X ( \beta^* - \beta^0) \|_2^2 + 2 \lambda \| \beta_{-S_0}^* \|_1 = 2 \lambda \beta_2^0 $$
       and with $C=4$ we get
     $$ \| X ( \beta^* - \beta^0) \|_2^2 + 2 \lambda \| \beta_{-S_0}^* \|_1 = 
    { \lambda \beta_2^0  }  + {  \lambda^2 \over \hat\varphi^2 }  . $$

     The next lemma has the situation of Lemma \ref{goodlasso2.lemma} but now with $C=1$
     instead of $C>1$.
     This  is an example where the minimizer of ${\cal L}  (\cdot)$  is not unique. 
     
      \begin{lemma}\label{goodlasso3.lemma}
    Let
    $$ X_3 = (X_1 + X_2 )/2 +V , \ X_4 = (X_1 + X_2)/2 -V , $$
    where $ V^T X_{S_0} =0 $.
        Then
    $$ \hat \phi^2 (S_0)=0 . $$
    Moreover when $\beta_1^0 \ge  \beta_2^0 \ge \lambda / \hat \varphi^2 $, we find
    that the vector
    $$\beta^* = \begin{pmatrix}  \beta_1^0 - \lambda / \hat \varphi^2 - \beta_3^* \cr
    \beta_2^0 - \lambda / \hat \varphi^2 - \beta_3^* \cr \beta_3^* \cr \beta_3^* 
     \end{pmatrix} $$
     is  for all $0 \le \beta_3^* \le \beta_2^0 - \lambda / \hat \varphi^2 $ a minimizer
     of ${\cal L}  (\cdot) $ and we have 
    $$ \| X ( \beta^* - \beta^0) \|_2^2 = {2 \lambda^2 \over \hat \varphi^2 }. $$
    $$ \| X ( \beta^* - \beta^0) \|_2^2 + 2 \lambda \| \beta_{-S_0}^* \|_1 \le 
    {4 \lambda \beta_2^0 }  - { 2 \lambda^2 \over \hat\varphi^2 } .  $$
     \end{lemma}

    \section{The case $s_0=m_0 =2N$} \label{s0=m0=2N.section}
  Suppose $S_0 = \{ 1 , \ldots , 2N \} $ and $S_0^c = \{2N+1 , \ldots , 4N \}$.
  We can easily extend the situation of Section \ref{s0=2=m0.section}, where $N=1$,
  to $N>1$ by assuming $N$ mutually orthogonal blocks of variables.
  This extension is trivial but nevertheless useful as it moves us away from
  a very low-dimensional situation.
  
  \begin{lemma}\label{s0=m0=2N.lemma}
  Set for $k=1 , \ldots , N$
  $$ ( X_{2k-1} , X_{2k} )^T ( X_{2k-1} , X_{2k} )= \begin{pmatrix} 1 & - \hat \rho_k \cr
  -\hat \rho_k & 1 \cr  \end{pmatrix} , \ \hat \varphi_k := 1- \hat \rho_k , $$
  and $(X_{2k-1} , X_{2k})$ orthogonal to $ \{ X_j \}_{j  \in S_0 \backslash \{ 2k-1, 2k \}} $. 
   Let for $k=1 , \ldots , N$
  $$X_{2N+2k-1}= C_k(X_{2k-1} + X_{2k} ) + U_k + V_k , 
  X_{2N+2k}= C_k(X_{2k-1} + X_{2k} ) + U_k - V_k , $$
  where $C_k >1$ and $ C_k^2 \hat \varphi^2/2 < 1$, $(U_k, V_k)$ orthogonal
  to $ \{ X_j \}_{j \in S_0 \backslash  \{ 2k-1, 2k \}} $ as well as to
  $\{ (U_j, V_j) \}_{j \not= k } $ , and $U_k^T V_k = 0$.
  Let $\hat \tau_k^2 := U_k^T U_k $ with $ 0< \hat \tau_k^2 < 1- C_k^2 \hat \varphi_k^2 $. 
  Then
 $$ \hat \Gamma^2 (S_0)= 2 \sum_{k=1}^N 1/ \hat \varphi_k^2  + \sum_{k=1}^N  C_k^2 / \hat \tau_k^2 .$$
 If, for $k=1 , \ldots , N$, $\beta_{2k-1}^0 \ge \beta_{2k}^0 \ge \lambda/ \hat \varphi^2 + \lambda
 C_k (C_k-1) / ( 2 \hat \tau_2)$ we obtain
 $$ \| X ( \beta^* - \beta^0) \|_2^2 + 2 \lambda \| \beta_{-S_0}^* \|_1 =
 \lambda^2 \hat \Gamma^2 (S_0) - \lambda^2 \| 1 / \hat \tau^2 \|_1 . $$
 
  \end{lemma}

\section{Further results with $s_0=2$}\label{Furthers0=2.section}
In the previous sections with $S_0 = \{ 1 , 2 \}$ we assume that each inactive variable is a given
a linear combination of the active ones plus an orthogonal term. In this section,
we assume the situation is the other way around: each active variable is a given
linear combination of the inactive ones plus an orthogonal term.

We first examine a case where the compatibility constant is zero,
  and the presence of non-active variables has big impact on the prediction
  error, even when the negative correlation $\hat \rho$ between active variables is small.
  Afterwards, this situation is slightly adjusted to one with positive compatibility constant,
  but the upper bounds are then a factor  too large.

 The next lemma has compatibility constant $\hat \phi^2 (S_0) $ equal to zero.

\begin{lemma} \label{m=s0=2.lemma} 
Let $S_0^c = \{ 3, 4 \} $ ($m_0=2$) and 
$$ (X_3 , X_4)^T (X_3,X_4)  = \begin{pmatrix} 1 & - \hat \theta \cr - \hat \theta & 1 \cr  \end{pmatrix}. $$
Assume that for some vector $(\gamma_3 , \gamma_4)^T= \gamma_{-S_0} \in \R^2$
with $1/2 < \gamma_3 < 1$ and $\gamma_4 = 1- \gamma_3$.
$$ X_1 = X_{-S_0} \gamma + V , \ X_2 = X_{-S_0} \gamma -V ,$$
where $X_{-S_0} := \{ X_j \}_{j \notin S_0}$ and where $V^T X_{-S_0} =0 $. Then
$$ \hat \phi^2 (S_0) = 0 $$
and
$$ \hat \varphi^2 = 2(1- 4 \gamma_3 (1- \gamma_3)) + 4 \gamma_3 (1- \gamma_3)\hat \psi^2 , $$
where $\hat \psi^2 := 1- \hat \theta $.\\
Furthermore,
 if $2\gamma_4  \beta_2^0 \ge \lambda / \hat \psi^2  $ we have
$$ \| X (\beta^* - \beta^0) \|_2^2 = {2 \lambda^2 \over \hat \psi^2}  $$
and
$$ \| X (\beta^* - \beta^0) \|_2^2 + 2 \lambda \| \beta_{-S_0}^* \|_1 =
4 \lambda \beta_2^0 - { 2 \lambda^2 \over \hat \psi^2 } \ge 4 \lambda \gamma_3 \beta_2^0 . $$
\end{lemma}

The above lemma illustrates that when the compatibility condition fails, the prediction error
$ \| X (\beta^* - \beta^0) \|_2^2$ can be as large as $4 \lambda \gamma_4 \beta_2^0$ 
where $\gamma_4 < 1/2$, even when the correlation $-\hat \rho$ between $X_1$ and $X_2$ is not close
to $-1$, i.e., even when $\hat \varphi^2$ is not close to zero
(as $\hat \varphi^2 > 2( 1- 4 \gamma_3 (1- \gamma_3)) $).

 We now consider two situations where the compatibility constant 
 is positive. Moreover, there are no false positives, i.e.\  $\| \beta_{-S_0}^* \|_1=0$.
 Indeed, in the two  Lemmas \ref{phi=rho.lemma} and 
 \ref{2comp>0.lemma}  the irrepresentable condition (\cite{ZY07}) holds.

 \begin{lemma}\label{phi=rho.lemma}
 Let $S_0^c: \{3 , 4 \}$ ($m_0=2$) and 
$$(X_3 , X_4)^T (X_3,X_4)   = \begin{pmatrix} 1 & - \hat \theta \cr - \hat \theta & 1 \cr  \end{pmatrix}$$
and write $\hat \psi^2 := 1 - \hat \theta$. Assume that
$$ X_1 = C(X_3 + X_4 )/2 + V , \ X_2 = C(X_3+X_4 )/2 -V ,$$
where $V^T X_{-S_0} =0 $ and $C>1$, $C^2 \hat \psi^2 / 2 < 1 $. 
Then 
$$\hat \phi^2 ( S_0) = (C-1)^2 \hat \psi^2 , \ \hat \Gamma^2 (S_0) = { 2 \over (C-1)^2 \hat \psi^2 } .$$
Moreover, $\hat \varphi^2  = C^2 \hat \psi^2$, and  
for $ \beta_2^0 \ge \lambda / \hat \varphi^2 $,
   $$ \| X (\beta^* - \beta^0) \|_2^2 = {2  \lambda^2  \over \hat \varphi^2  } =
    \lambda^2 \hat \Gamma^2 (S_0) { (C-1)^2 \over C^2 }
 $$
   and $\| \beta_{-S_0}^* \|_1 =0$. 
\end{lemma}

In other words, the upper bound $\lambda^2 |S_0| / \hat \phi^2 (S_0)$ is a factor $C^2/(C-1)^2$ too large
   in this case.
   
   In the last result of this paper, we again let $s_0=2$ but now $m_0$ is arbitrary. Moreover, we assume
   that the inactive variables are orthogonal to each other.

    \begin{lemma} \label{2comp>0.lemma} Let 
 $S_0 = \{1 , 2 \}$,  $\hat \Sigma_{-S_0 , -S_0 } =I$ and 
 $$X_1= C X_{-S_0} \gamma_{-S_0} + V, \ X_2 = X_{-S_0} \gamma_{-S_0} -V, $$
  where $X_{-S_0} := \{ X_j \}_{j \notin S_0} $ and where 
 $V^T X_{-S_0} = 0$. Assume moreover
 $\| \gamma_{-S_0} \|_1 = 1 $, $2 C^2 \| \gamma_{-S_0 } \|_2^2 < 1$ and $\| \gamma_{-S_0} \|_{\infty} \le C
 \| \gamma_{-S_0} \|_2^2 $. 
  Then 
 $$ \hat \phi^2 (S_0)=  2 \min_{\| \beta_{-S_0}\|_1 \le 1 } \| C \gamma_{-S_0} - \beta_{-S_0} \|_2^2,
 \ \hat \varphi^2 = 2 C^2 \| \gamma_{-S_0} \|_2^2 , $$
 and moreover for $\beta_1^0 \ge \beta_2^0 \ge \lambda / \hat \varphi^2 $
 $$ \| X ( \beta^* - \beta^0 ) \|_2^2  = {2 \lambda^2  \over \hat \varphi^2 }.    $$
  \end{lemma}
 
      \hfill $\sqcup \mkern -12mu \sqcap$
    
    \begin{corollary}\label{2comp>0.corollary}
    An example of a vector $\gamma_{-S_0}$ and constant $C$ in Lemma \ref{2comp>0.lemma} is
$$ \gamma_{-S_0} = ( \underbrace{1 /m_0 , \cdots ,1 /m_0}_{m_0 \ \times}  )^T ,$$
and $1 <C^2 < m_0/2    $. 
Then
$$  \hat \phi^2 (S_0)= 2 (C-1)^2 / m_0, \ \hat \Gamma^2 (S_0)= { m_0 \over (C-1)^2 }$$  and
  $$ \| X ( \beta^* - \beta^0 ) \|_2^2  = { \lambda^2  m_0 \over C^2  } = 
  \lambda^2 \hat \Gamma^2 (S_0)  { (C-1)^2 \over C^2 }  .  $$
  So again there is a gap with Lemma \ref{upperbound1.lemma}, but it is small for $C$ large. 
\end{corollary}

  \section{Proofs} \label{proofs.section}
  
  In the proofs, we sometimes use the following notation. 
 The matrix with columns in $S\subset \{ 1 , \ldots , p\}$ is written as $X_S := \{ X_j \}_{j \in S} $
and $X_{-S} := \{ X_j \}_{j \notin S} $ has its columns in $S^c$.
The order in the columns is taken increasing in the index (i.e., we remove some columns and otherwise keep the original
ordering). 
We write
$$\hat \Sigma_{S,S} := X_S^T X_S   , \ \hat \Sigma_{-S,S}: = X_{-S}^T X_S  , $$ $$
\hat \Sigma_{S,-S} := X_S^T X_{-S}  , \ \hat \Sigma_{-S,-S} := X_{-S}^T X_{-S}   . $$

  In the proofs of results from Section \ref{p=s0=2.section} and onwards we present explicit expressions for the minimizer $\beta^*$ showing it
  is the solution of the KKT conditions. 
 One may check that  the solution is unique in each case
except for Lemma \ref{goodlasso3.lemma}.

 \subsection{Proof of the results in Section \ref{noisy.section}}
  
  Theorem \ref{silly.theorem} and its proof are stated as Problem 2.4 in \cite {vdG2016}. Here, we present a
  complete proof. 
  For this we need some auxiliary lemmas.

 \begin{lemma}\label{hidden-oracle1.lemma} It holds that
 $$ \| X({\hat \beta} - {\beta^*} )\|_2^2 +  \lambda \| \hat \beta \|_1 -
  \lambda \hat \beta^T z^* \le   (\hat \beta -\beta^* )^T X^T   \epsilon.$$
    \end{lemma}
    
    {\bf Proof of Lemma \ref{hidden-oracle1.lemma}.}
 By the KKT conditions for $\hat \beta$
 $$ - X^T (Y- X \hat \beta)  + \lambda \hat z =0 , $$
 where $\hat z \in \partial \| \hat \beta \|_1 $.
 In other words
 $$ \hat \Sigma ( \hat \beta - \beta^0 ) + \lambda \hat z = X^T \epsilon  . $$
 By the KKT conditions for $\beta^*$
 $$ \hat \Sigma  ( \beta^* - \beta^0 ) + \lambda z^* =0 . $$
 Hence, taking the difference
 $$ \hat \Sigma ( \hat \beta - \beta^* ) + \lambda (\hat z -  z^* ) = X^T \epsilon  . $$
 Multiply by $(\hat \beta - \beta^*)^T$ to find
 $$\| X({\hat \beta } - {\beta^*} )\|_2^2 + \lambda ( \hat \beta - \beta^* )^T (\hat z - z^* ) =
  (\hat \beta -\beta^* )^T X^T   \epsilon . $$
 But
 $$( \hat \beta - \beta^* )^T (\hat z - z^* )= \hat \beta^T (\hat z - z^*) + \beta^{*T} (z^* - \hat z) . $$
 Both terms are non-negative: since $\hat \beta^T z^* \le 
 \| \hat \beta \|_1 \| z^* \|_{\infty} \le \| \hat \beta \|_1$ we have 
 $$\hat \beta^T (\hat z - z^*)= \| \hat \beta \|_1 - \hat \beta^T z^* \ge 0 $$
 and by the same argument
 $$ \beta^{*T} (z^* - \hat z) = \| \beta^* \|_1 - \beta^{*T} \hat z \ge 0 . $$
 Dropping the term $\| \beta^* \|_1 - \beta^{*T} \hat z$ therefore yields 
 $$ \| X({\hat \beta} - {\beta^*} )\|_2^2 +  \lambda \| \hat \beta \|_1 -
  \lambda \hat \beta^T z^* \le  (\hat \beta -\beta^* )^T X^T   \epsilon . $$
  \hfill $\sqcup \mkern -12mu \sqcap$
  
  Recall that the vector $\beta^*$ satisfies the KKT conditions
 $$ \hat \Sigma ( \beta^* - \beta^0)  + \lambda z^* =0 , $$
 where $z^* \in \partial \| \beta^* \|_1 $. 
 
 Define 
 $$ \bar S_*:= \{ j:\ |z_j^* | \ge 1-\eta \} . $$
 Note that $\bar S_* \supset S_*$ where $S_* $ is the active set of $\beta^*$.
 We write $\bar s_* := | \bar S_* | $.
 
 \begin{lemma} \label{sizeS*.lemma}
 It holds that
 $$\bar s_* \le  {  \Lambda_{\rm max}(\hat \Sigma)   \over \lambda^2 (1- \eta)^2 }  \| X ({\beta^*} - \beta^0)  \|_2^2. $$
  \end{lemma}
 
 {\bf Proof of Lemma \ref{sizeS*.lemma}.} 
 By the KKT conditions for $\beta^*$ it is true that
 $$ \lambda^2 \|  z^* \|_2^2 = (\beta^* - \beta^0)^T \hat \Sigma^2 (\beta^* - \beta^0) \le
    \Lambda_{\rm max} (\hat \Sigma) \|X({\beta^*} - \beta^0) \|_2^2 . $$
 On the other hand
 $$ \| z^* \|_2 ^2 \ge \|  z_{\bar S_*}^* \|_2^2 \ge (1- \eta)^2 \bar s_* . $$
 Hence
 $$ \bar s_*  \le{ \Lambda_{\rm max} (\hat \Sigma)    \over  \lambda^2 (1-\eta)^2 }
 \| X( \beta^* - \beta^0)  \|_2^2  . $$
 \hfill $\sqcup \mkern -12mu \sqcap$

 Define the random variable 
 $$  {\bf V}^2 (\bar S_*) := \max_{ \| X{\beta_{\bar S_*} } \|_2 =1 } | \beta_{\bar S_*}^TX^T \epsilon   |
     .$$
     Define moreover the vector $X_{\bar S_*} \hat \gamma_{\bar S_*} $ as the projection of
 $X \hat \beta $ on the space spanned by the columns of $X_{\bar S_*} $ and let ${\bf w}$ be the random variable
 $$ {\bf w} := {
 \biggl  \| \biggl [X\hat \beta- X_{\bar S_*} \hat \gamma_{\bar S_*} \biggr ]^T  \epsilon \biggr \|_{\infty} /
\  \| \hat \beta_{-\bar S_*} \|_1 }.  $$

   \begin{lemma}\label{hidden-oracle2.lemma} We have
  $$\| X({\hat \beta} - {\beta^* } )\|_2^2 +  2 (\eta \lambda - {\bf w} ) \| \hat \beta_{-\bar S^* } \|_1  \le
 {\bf V}^2 (\bar S_*) . $$
   \end{lemma}

 {\bf Proof of Lemma \ref{hidden-oracle2.lemma}.}
    By Pythagoras' theorem, and using that $S_* \subset \bar S_*$
 $$ \| X(\hat \beta -\beta^* )\|_2^2 = \| X_{\bar S_*} \hat \gamma_{\bar S_*} - X \beta^* \|_2^2 +
 \| X\hat \beta - X_{\bar S_*} \hat \gamma_{\bar S_*} \|_2^2 . $$
 Therefore, in view of Lemma \ref{hidden-oracle2.lemma},
 $$  \| X_{\bar S_*} \hat \gamma_{\bar S_*} - X \beta^* \|_2^2 +
 \| X \hat \beta - X_{\bar S_*} \hat \gamma_{\bar S_*} \|_2^2 
 +  \lambda \| \hat \beta \|_1 -
  \lambda \beta^T z^* $$ $$ \le  \biggl [X_{\bar S_*} \hat \gamma_{\bar S_*} - X \beta^*  \biggr ]^T \epsilon+
   \biggl [ X\hat \beta - X_{\bar S_*} \hat \gamma_{\bar S_*} \biggr ]^T  \epsilon . $$
 By the Cauchy-Schwarz inequality
 $$\biggl [X_{\bar S_*} \hat \gamma_{\bar S_*} - X \beta^*  \biggr ]^T \epsilon
  \le {\bf V} (\bar S_*)  \| X_{\bar S_*} \hat \gamma_{\bar S_*} - X \beta^* \|_2. $$
 Moreover, by the definition of ${\bf w}$
 $$  \biggl [ X\hat \beta - X_{\bar S_*} \hat \gamma_{\bar S_*} \biggr ]^T  \epsilon \le {\bf w}   \| \hat \beta_{-\bar S_*} \|_1 . $$
 On the other hand, $|z_j^* | \le 1- \eta $ for all $j \notin \bar S_* $ and hence
 $$ \| \hat \beta_{-\bar S_*} \|_1 - z_{-\bar S_*}^{*T} \hat \beta_{-\bar S_*} \ge \eta 
 \| \hat \beta_{-\bar S_*} \|_1 . $$
 We thus arrive at
 \begin{eqnarray*}
\| X_{\bar S_*} \hat \gamma_{\bar S_*} - X \beta^* \|_2^2 & +& 
 \| X\hat \beta - X_{\bar S_*} \hat \gamma_{\bar S_*} \|_2^2 
 +  \eta \lambda \| \hat \beta_{-\bar S^*}  \|_1 \\
 &\le&   {\bf V} (\bar S_*)  \| X_{\bar S_*} \hat \gamma_{\bar S_*} - X \beta^* \|_2 + 
  {\bf w}  \| \hat \beta_{-\bar S_*} \|_1  \\
  & \le & {\bf V}^2 (\bar S_*) / 2 + \| X_{\bar S_*} \hat \gamma_{\bar S_*} - X \beta^* \|_2^2/2 +
   {\bf w}  \| \hat \beta_{-\bar S_*} \|_1 
   \end{eqnarray*} 
  or
 $$ \| X_{\bar S_*} \hat \gamma_{\bar S_*} - X \beta^* \|_2^2 + 2\| X\hat \beta - X_{\bar S_*} \hat \gamma_{\bar S_*} \|_2^2
  + 2(\eta \lambda - {\bf w}) \| \hat \beta_{-\bar S^* } \|_1 
 \le {\bf V}^2 (\bar S_*)   . $$
 But then also
 $$\| X ( \hat \beta -\beta^*)  \|_2^2 +  2 (\eta \lambda - {\bf w} ) \| \hat \beta_{-\bar S^* } \|_1  \le
 {\bf V}^2 (\bar S_*) . $$
 
 \hfill $\sqcup \mkern -12mu \sqcap$

 \begin{lemma} \label{project-dual-norm.lemma} 
 Let $ \lambda_0 := \sqrt {2 \log(2p/\alpha) / n } $.
 It holds that with probability at least $1- \alpha$ that 
 $$ {\bf w} \le \lambda_0. $$
\end {lemma}

{\bf Proof of Lemma \ref{project-dual-norm.lemma}.}
 Write the singular value decomposition of  $X_{\bar S_*} $ as 
 $$X_{\bar S_*} = \bar P_* \bar \Lambda_*^{1/2} \bar Q_* $$
 where
 $ \bar P_*^T\bar P_* =  I $, $\bar Q_*^T \bar Q_* = I $ and
 $\bar \Lambda_* $ the diagonal matrix of eigenvalues of $X_{\bar S_*}^T X_{\bar S_*} $.
 Since $X\hat \beta = X_{\bar S_*} \hat \beta_{\bar S_*} + X_{-\bar S_*} \hat \beta_{-\bar S_*} $, we see that
 $$X\hat \beta - X_{\bar S_*} \hat \gamma_{\bar S_*}  = (I- \bar P_* \bar P_*^T ) X_{-\bar S_*} \hat \beta_{-\bar S_*} .$$
 Hence
 $$\biggl  [X\hat \beta- X_{\bar S_*} \hat \gamma_{\bar S_*} \biggr ]^T  \epsilon =
( X_{-\bar S_*} \hat \beta_{-\bar S_*} )^T  ( I- \bar P_* \bar P_*^T )  \epsilon . $$
Thus
$$ {\bf w} \le  \| X_{-\bar S_*}^T  ( I- \bar P_* \bar P_*^T )  \epsilon\|_{\infty} . $$
 The diagonal elements of the matrix
 $$ X_{-\bar S_*}^T (I- \bar P_* \bar P_*^T)X_{-{\bar S_*} }  $$
 are projected versions of the columns of
 $X_{-{\bar S_*} } $ and hence at most $\max_{j \in \bar S_*} \| { X }_j \|_2^2  $, which is by 
 assumption at most $1$. 
  It follows that each element of the vector  $ \sqrt n X_{-\bar S_*}^T   ( I- \bar P_* \bar P_*^T )  \epsilon $
 is
 normally distributed with mean zero and variance at most $1$. The dimension of this vector is
 at most $p$. Now use that for standard normal random variables $W_1 , \ldots , W_p$, and for any $t>0$,
 \begin{eqnarray*}
   \PP ( \max_{1 \le j \le p} | W_j | > \sqrt { 2 (\log (2p) + t) } )  
 &\le&
 p \PP(  | W_1 | \ge \sqrt { 2 (\log (2p) + t) } )  \\
 & \le&  2 p \exp[- (\log (2p +t)] = \exp[-t] . 
 \end{eqnarray*}
 Apply this with $t = \log (1/ \alpha) $. 
 \hfill $\sqcup \mkern -12mu \sqcap$

  \begin{lemma} \label{chi-square2.lemma} We have 
   $$ \PP ({\bf V} (\bar S_*)  \ge \sqrt {\bar s_* /n} + \sqrt {2\log (1/\alpha_1) /n} ) \le \alpha_1  . $$
 \end{lemma}

 {\bf Proof of Lemma \ref{chi-square2.lemma}.}
 Let $\chi_T^2$ be  chi-squared random variable with $T$ degrees of freedom.
 Lemma 1 in \cite{laurent2000adaptive} says that
for all $t>0$
 $$ \PP ( \chi_T^2 \ge T + 2 \sqrt { Tt} + 2t ) \le \exp[-t] . $$
 Since $T + 2 \sqrt {Tt} + 2t  \le (\sqrt {T} + \sqrt {2t} )^2 $
 we find
 $$ \PP (\chi_T \ge \sqrt {T} + \sqrt {2t} ) \le \exp[-t] . $$
 Apply this with $t = \log (1/ \alpha_1) $. 
   \hfill $\sqcup \mkern -12mu \sqcap$

 {\bf Proof of Theorem \ref{silly.theorem}.} 
 We know by Lemma \ref{project-dual-norm.lemma} that with probability at least $1- \alpha$
 $$ {\bf w} \le \lambda_0 $$
 and from 
 Lemma \ref{chi-square2.lemma}, with probability at least $1- \alpha_1$
 $$ {\bf V} (\bar S_*) \le  \sqrt {  \bar s_* \over n } + \sqrt{2 \log (1/\alpha_1) \over n} $$
 By Lemma \ref{sizeS*.lemma}
 $$ \bar s_*  \le { \Lambda_{\rm max} (\hat \Sigma)    \over  \lambda^2 (1- \eta)^2 } 
 \| X( \beta^* - \beta^0)  \|_2^2 . $$
 Hence with probability at least $1- \alpha_1$
 $${\bf V}(\bar S_*) \le \sqrt {  \Lambda_{\rm max} (\hat \Sigma) \over  n\lambda^2 (1- \eta)} \| X({\beta^*} - {\beta^0}) \|_2 +
 \sqrt {2 \log (1/\alpha_1) \over n  } . $$
 Combine this with Lemmas \ref{hidden-oracle1.lemma}  and \ref{hidden-oracle2.lemma} and invoke
 the condition $\eta \lambda > \lambda_0$ to complete the proof.
 \hfill $ \sqcup \mkern -12mu \sqcap$

  \begin{lemma} \label{sizeS*bis.lemma}
 Suppose that
 $$
 \| \hat \Sigma - \Sigma_0\|_{\infty}  \| \beta^* - \beta^0 \|_1< \lambda (1- \eta) .$$
 Then
 $$\bar s_* \le  { \Lambda_{\rm max} ( \Sigma_0 )  \biggl ( \| X( \beta^* - \beta^0)  \|_2^2  +  \| \hat \Sigma - \Sigma_0\|_{\infty} \| \beta^* - \beta^0 \|_1^2 \biggr ) \over 
  \biggl (  \lambda (1- \eta )- \| \hat \Sigma - \Sigma_0\|_{\infty}  \| \beta^* - \beta^0 \|_1 \biggr )^2 }  .  $$
  \end{lemma}
 
 {\bf Proof of Lemma \ref{sizeS*bis.lemma}.} 
 We start again with the KKT conditions for $\beta^*$ 
  $$ \hat \Sigma  ( \beta^* - \beta^0 ) + \lambda z^* =0 . $$
  Then
  $$ \Sigma_0 ( \beta^*- \beta^0) + ( \hat \Sigma - \Sigma_0) (\beta^* - \beta^0) =
  -\lambda z^* . $$
  But for all $j$
  $$ | (( \hat \Sigma - \Sigma_0) (\beta^* - \beta^0))_j| \le \| \hat \Sigma - \Sigma_0\|_{\infty}  \| \beta^* - \beta^0 \|_1 $$
    so 
  $$|  ( \Sigma_0 ( \beta^*- \beta^0))_j | \ge \lambda | z_j^*| - 
  \| \hat \Sigma - \Sigma_0\|_{\infty}  \| \beta^* - \beta^0 \|_1 .$$
  If $| z_j^* | > (1- \eta)$ we get 
  $$\lambda | z_j^*| - 
  \| \hat \Sigma - \Sigma_0\|_{\infty}  \| \beta^* - \beta^0 \|_1 > \lambda (1- \eta )- \| \hat \Sigma - \Sigma_0\|_{\infty}  \| \beta^* - \beta^0 \|_1 > 0 . $$
  Thus
  $$ \sum_{j \in {\bar S_*} } \biggl (  \lambda z_j -\| \hat \Sigma - \Sigma_0\|_{\infty}  \| \beta^* - \beta^0 \|_1 \biggr )^2
  \ge \bar s_* \biggl (  \lambda (1- \eta )- \| \hat \Sigma - \Sigma_0\|_{\infty}  \| \beta^* - \beta^0 \|_1 \biggr )^2 . $$
  On the other hand
 $$ \sum_{j \in \bar S_*} | (\Sigma_0 ( \beta^* - \beta^0))_j  |^2 \le 
 \Lambda_{\rm max} ( \Sigma_0) (\beta^*- \beta^0)^T \Sigma_0 ( \beta^* - \beta^0)  . $$
 $$ \le  \Lambda_{\rm max} ( \Sigma_0) \biggl ( \| X(\beta^* - \beta^0) \|_2^2 +\| \hat \Sigma - \Sigma_0 \|_{\infty} 
 \| \beta^* - \beta^0 \|_1^2 \biggr ) .$$
  Hence
 $$ \bar s_*  \le { \Lambda_{\rm max} ( \Sigma_0 )  \biggl ( \| X( \beta^* - \beta^0)  \|_2^2  + 
 \| \hat \Sigma - \Sigma_0 \|_{\infty}  \| \beta^* - \beta^0 \|_1^2 \biggr ) \over 
  \biggl (  \lambda (1- \eta )- \| \hat \Sigma - \Sigma_0\|_{\infty}  \| \beta^* - \beta^0 \|_1 \biggr )^2 }  
 . $$
 \hfill $\sqcup \mkern -12mu \sqcap$
 
 {\bf Proof of Theorem \ref{sillybis.theorem}.}
 We have by Lemma \ref{sizeS*bis.lemma}
$$ \bar s_*\le  { \Lambda_{\rm max} ( \Sigma_0 )  \biggl ( \| X( \beta^* - \beta^0)  \|_2^2  + 
\| \hat \Sigma - \Sigma_0 \|_{\infty}   \| \beta^* - \beta^0 \|_1^2 \biggr ) \over 
  \biggl (  \lambda (1- \eta )- \| \hat \Sigma - \Sigma_0\|_{\infty}  \| \beta^* - \beta^0 \|_1 \biggr )^2 }  . $$
  So with probability at least $1- \alpha_1$,
  \begin{eqnarray*}
 {\bf V}(\bar S_*) &\le &{ \Lambda_{\rm max}^{1/2}  ( \Sigma_0 )  \biggl ( \| X( \beta^* - \beta^0)  \|_2^2  + 
  \| \hat \Sigma - \Sigma_0 \|_{\infty}  \| \beta^* - \beta^0 \|_1^2 \biggr )^{1/2}  \over 
  \biggl (  \lambda (1- \eta )- \| \hat \Sigma - \Sigma_0\|_{\infty}  \| \beta^* - \beta^0 \|_1 \biggr )} \\
  & + &
 \sqrt {2 \log (1/\alpha_1) \over n  } . 
 \end{eqnarray*} 
 The proof can be completed along the same lines as the proof of Theorem \ref{silly.theorem}.
 \hfill $\sqcup \mkern -12mu \sqcap$

 \subsection{Proofs for Section \ref{upper-bound.section}}
     
  {\bf Proof of Lemma \ref{upperbound1.lemma}.}
 By the KKT conditions
 $$ \hat \Sigma ( \beta^* - \beta^0) + \lambda z^* = 0 , \ z^* \in \partial \| \beta^* \|_1 .$$
 Hence
  \begin{eqnarray*}\label{basic.equation}
 0 \le  \| X ( \beta^* - \beta^0) \|_2^2  &= &
 (\beta^* - \beta^0)^T  \hat \Sigma ( \beta^* - \beta^0)\\
 &=& \lambda ( \beta^0 - \beta^* )^T z^* \\
 & \le & \lambda \| \beta^0 \|_1 - \lambda \| \beta^* \|_1 .
 \end{eqnarray*}
 Therefore the first bound of the lemma holds.
 Continuing with (\ref{basic.equation}) and applying the definition of the compatibility
 constant $\hat \phi^2 (S_0)$ one finds
  \begin{eqnarray*}\
 0 \le  \| X ( \beta^* - \beta^0) \|_2^2  &\le &
 \lambda \| \beta^0 \|_1 - \lambda \| \beta^* \|_1 \\
 & \le & \lambda \| \beta^0 - \beta_{S_0}^* \|_1 - \lambda \| \beta_{-S_0}^* \|_1 \\
 & \le & \lambda  \sqrt {s_0}  \| X ( \beta^* - \beta^0) \|_2   \hat \phi (S_0) -  \lambda \| \beta_{-S_0}^* \|_1\\
 & \le & {\lambda^2 s_0 /( 2 \hat \phi^2 (S_0) )} + \| X ( \beta^* - \beta^0) \|_2^2 /2 
 -  \lambda \| \beta_{-S_0}^* \|_1 .
 \end{eqnarray*}
 This yields the second bound of the lemma.
 \hfill $\sqcup \mkern -12mu \sqcap$

 {\bf Proof of Lemma \ref{upperbound.lemma}.}
  The first minimum $\hat {\cal U}_I (\beta^0)$ for the prediction error follows from Lemma
  \ref{upperbound1.lemma}.
  
  We recall the KKT conditions
  $$ \hat \Sigma ( \beta^* - \beta^0) + \lambda z^* = 0 , \ z^* \in \partial \| \beta^* \|_1 .$$
  For the second minimum $\hat {\cal U}_{II} (\beta^0) $, let $S\subset \{ 1 , \ldots , p \}$ be arbitrary. We note that
 when $\| X ( \beta^* - \beta^0 )\|_2^2 - 2 \lambda \| \beta_{-S}^0 \|_1 \le 0 $ there is
 nothing to prove here. So let us assume $\| X ( \beta^* - \beta^0) \|_2^2 - 2 \lambda \| \beta_{-S}^0 \|_1 \ge 0 $.
 Then we have by the KKT conditions
 \begin{eqnarray*}
 0 &\le& \| X ( \beta^* - \beta^0 )\|_2^2 - 2 \lambda \| \beta_{-S}^0 \|_1 \\
& \le &\lambda  \| \beta^0\|_1- \lambda \| \beta^* \|_1 - 2 \lambda \| \beta_{-S}^0 \|_1 \\
& \le & \lambda \| \beta_S^* - \beta_S^0 \|_1 - \lambda \| \beta_{-S}^* \|_1 - \lambda \| \beta_{-S}^0 \|_1 \\
& \le & \lambda \| \beta_S^* - \beta_S^0 \|_1 - \lambda \| \beta_{-S}^* - \beta_{-S}^0 \|_1 .
 \end{eqnarray*}
 By the definition of the compatibility constant we now find
 \begin{eqnarray*} 
& \ &  \| X ( \beta^* - \beta^0 ) \|_2^2 - 2 \lambda \| \beta_{-S}^0 \|_1\\
 &  \le &
 \lambda \sqrt {|S|}\| X (\beta^* - \beta^0) \|_2  /  \hat \phi (S)   - \lambda \| \beta_{-S}^* \|_1
 - \lambda \| \beta_{-S}^0 \|_1\\
 &\le & \lambda \sqrt {|S|} \| X (\beta^* - \beta^0) \|_2 / \hat \phi (S) - \lambda \| \beta_{-S}^0 \|_1 .
 \end{eqnarray*}
 It follows that
 $$ \biggl ( \| X ( \beta^* - \beta^0 ) \|_2  - \lambda \sqrt |S| / ( 2 \hat \phi ( S)) \biggr )^2 $$
 $$ \le \lambda^2 |S| / (4 \hat \phi^2 ( S)) + \lambda \| \beta_{-S} \|_1 . $$

 We now turn to the third minimum $\hat {\cal U}_{III} (\beta^0) $. For any $\beta$ 
 $$ ( \beta^* - \beta )^T \hat \Sigma ( \beta^* - \beta^0) + ( \beta^* - \beta )^T z^* = 0 . $$
 We have
 \begin{eqnarray*}
 ( \beta^* - \beta )^T \hat \Sigma ( \beta^* - \beta^0)=
 &\|  X( \beta^* - \beta^0 )\|_2^2 / 2 &- \| X( \beta - \beta^0 )\|_2^2 / 2 \\ 
 + & \|  X( \beta^* - \beta)  \|_2^2 / 2 &. 
 \end{eqnarray*}
 Let $S \subset \{ 1 , \ldots , p \}$.
 If  $( \beta^* - \beta )^T \hat \Sigma ( \beta^* - \beta^0) - 2 \lambda \| \beta_{-S} \|_1 \le 0 $
 we are done.
 On the other hand, if
 $( \beta^* - \beta )^T \hat \Sigma ( \beta^* - \beta^0) - 2 \lambda \| \beta_{-S} \|_1 \ge 0 $
 we get
\begin{eqnarray*}
0 & \le & ( \beta^* - \beta )^T \hat \Sigma ( \beta^* - \beta^0) - 2 \lambda \| \beta_{-S} \|_1 \\
 &\le &  \lambda \| \beta_S \|_1 - \lambda \| \beta^* \|_1 - \lambda \| \beta_{-S} \|_1 \\
 & \le &\lambda \| \beta_{S}^* - \beta_S \|_1 - \lambda \| \beta_{-S}^* \|_1- \lambda \| \beta_{-S} \|_1 \\
 & \le &\lambda \| \beta_{S}^* - \beta_S \|_1 - \lambda \| \beta_{-S}^*- \beta_{-S} \|_1 .
 \end{eqnarray*}
 We can apply the definition of the compatibility constant to find
 \begin{eqnarray*}
 ( \beta^* - \beta )^T \hat \Sigma ( \beta^* - \beta^0) 
& \le & \lambda \sqrt { |S| } \| X( \beta^* - \beta) \|_2/  \hat \phi (S) + \lambda \| \beta_{-S} \|_1 \\
& \le & \lambda^2 |S| / (2 \hat \phi^2 (S)  ) + \|  X( \beta^* - \beta) \|_2^2 / 2 +\lambda \| \beta_{-S} \|_1 ,
\end{eqnarray*}
which gives
$$ 
\| X( \beta^* - \beta^0 )\|_2^2 / 2 - \| X( \beta - \beta^0 )\|_2^2 / 2 +
 \| X( \beta^* - \beta)  \|_2^2 / 2  $$
 $$
  \le  \lambda^2 |S| / (2 \hat \phi^2 (S)  ) + \| X( \beta^* - \beta) \|_2^2 / 2  + \lambda \| \beta_{-S} \|_1. 
$$

\hfill $\sqcup \mkern -12mu \sqcap$ 

\subsection{Proof of the lemma in Section \ref{design.section}}

 {\bf Proof of Lemma \ref{1fair.lemma}.}
   Suppose on the contrary that $\hat \phi (\{1\})=0$. Then there exists a $ \gamma_{-1}$ with
   $\| \gamma_{-1}  \|_1 = 1$ such that $X_1= X_{-1} \gamma_{-1} $.
   This gives 
   $$ 1= \| X_{1} \|_2^2 = \| X_{-1} \gamma_{-1} \|_2^2   .$$
   We show that this is not possible.
   We let $X_{-1}$ be an $n \times m_0$-matrix
and prove the result by induction in $m_0$.\\
$\circ$ $m_0=1$: Trivial.\\
$\circ$ $m_0=2$: Let $ \hat \vartheta := X_2^T X_3$. Assume without loss of generality that
$\gamma_{-1}^T= (\gamma_2 , \gamma_3)$ has both its components non-negative. Then
$\gamma_3 = 1- \gamma_2$ and 
\begin{eqnarray*}
\| X_{-1} \gamma_{-1} \|_2^2 & =& \gamma_2^2 + (1- \gamma_2 )^2 + 2 \gamma_2 (1-\gamma_2) \hat \vartheta\\
&=& 1+ 2 \gamma_2 (1- \gamma_2) \hat \vartheta
\end{eqnarray*}
This can only be equal to 1 if $\gamma_2=0$ or $\gamma_2= 1$ or $ \hat \vartheta=-1$, all cases
which we excluded.\\
$\circ$ Induction step: suppose it is true for the value $m_0-1$:
for all $\tilde \gamma_{-1} $ with $\tilde \gamma_{j_0 } = 0$ for some
$j_0 \in \{ 2 , \ldots , m_0+1 \} $ and with $\| \tilde \gamma  \|_1 =1$ it holds that
$ \| X_{-1} \tilde \gamma_{-1} \|_2^2   < 1 $.
 Let $\gamma_{-1}^T =
(\gamma_2 , \ldots , \gamma_{m_0+1})$ be a vector with
$\| \gamma_{-1} \|_1 = 1 $ and with $|\gamma_{m_0+1}| < 1$. Then we know
by induction that either $\| X_{-1} \gamma_{-1} - X_{m_0+1} \gamma_{m_0+1}  \|_2/  <  (1- |\gamma_{m_0+1} | )  $ or
there is a $j_0 \in \{ 2 , \ldots , m_0 \}$ such that $|\gamma_{j_0} |= 1- |\gamma_{m_0+1}|$. In the last case
all values $j \in \{ 2 , \ldots , m_0 \} $ other than $ j_0$ must be zero so it brings us back to the case $m_0=2$. In the first case
we have by the triangle inequality
\begin{eqnarray*}
\| X_{-1} \gamma_{-1} \|_2  &\le& \| X_{-1} \gamma_{-1} - X_{m_0+1} \gamma_{m_0+1}  \|_2 +  |\gamma_{m_0+1}| \\
&<& (1- |\gamma_{m_0+1}|  ) + |\gamma_{m_0+1} |=1 .
\end{eqnarray*}
\hfill $\sqcup \mkern -12mu \sqcap$

\subsection{Proofs for Section \ref{p=s0=2.section}}

{\bf Proof of Lemma \ref {hatphi1.lemma}.} 
The coefficient of the projection of $X_1$ on $X_2$ is
$$ \arg \min_{\beta_2 \in \R}  \| X_1 - X_2 \beta_2 \|_2^2   = - \hat \rho . $$
Since $0 < \hat \rho  < 1 $ we thus find
$$ \hat \phi^2 ( \{ 1 \} ) := \min_{| \beta_2 | \le 1 }  \| X_1- X_2 \beta^2 \|_2^2 =
\| X_1 + \hat \rho X_2 \|_2^2 = 1- \hat \rho^2 . $$
As $\hat \varphi^2 = (1- \hat \rho)$ we have 
\begin{eqnarray*}
1- \hat \rho^2 &=& (1- \hat \rho) (1+ \hat \rho) \\
&=&  (1- \hat \rho) (2- (1- \hat \rho) ) \\
&=&  \hat \varphi^2 (2- \hat \varphi^2  ). 
\end{eqnarray*}
 The second result follows from symmetry arguments:
 the minimum of $\| X_1 \beta_1 + X_2 \beta_2$ over $|\beta_1 | + | \beta_2 |=1$ is reached at
 equal values for $\beta_1$ and $\beta_2$. 
\hfill $\sqcup \mkern -12mu \sqcap$

{\bf Proof of Lemma \ref{p=2-exact.lemma}.}
One readily verifies that $0 \le \beta_1^* \le \beta_1^0$ and $0 \le \beta_2^* \le \beta_2^0$. 
Let $\Delta_1:= \beta_1^0 - \beta_1^* $ and $\Delta_2:= \beta_2^0 - \beta_2^* $.
 Recall the KKT conditions
  $$  \hat \Sigma \Delta =  \lambda z^*  , \ z^* \in \partial \| \beta^* \|_1 .$$
  $\circ$  ${\rm Case \ 1  } $:
 \fbox{$\lambda/\hat \varphi^2  \le   \beta_2^0 $}. 
For $\Delta_1=\Delta_1= \lambda/ \hat \varphi^2$
$$\hat \Sigma \Delta = \lambda \begin{pmatrix} 1 \cr 1 \cr \end{pmatrix} . $$
Since
 $$ \begin{pmatrix} \beta_1^* \cr \beta_2^* \cr \end{pmatrix} = 
 \begin{pmatrix} \beta_1^0 - \Delta_1  \cr \beta_2^0 - \Delta_2  \cr \end{pmatrix} $$
 has both its components non-negative, it 
 is a solution of the KKT conditions, in fact it is the unique solution.
 
 $\circ$ ${\rm Case \ 2 } $:  \fbox{$ \beta_2^0 <  \lambda/\hat \varphi^2  \le  \beta_2^0 + (\beta_1^0 - \beta_2^0 )/
\hat \varphi^2  $}.
With
$\Delta_1= \lambda+ (1 - \hat \varphi^2 )  \beta_2^0  = \lambda + \hat \rho \beta_2^0 $ and $\Delta_2=  \beta_2^0 $ we obtain
$$ \hat \Sigma \Delta = \begin{pmatrix} 1 & -\hat \rho\cr - \hat \rho & 1 \cr \end{pmatrix} \begin{pmatrix}
\lambda+ \hat \rho  \beta_2^0 \cr  \beta_2^0 \cr  \end{pmatrix}= 
\lambda \begin{pmatrix}  1 \cr z_2^* \end{pmatrix} $$
with $z_2^*=  - \hat \rho + (1 - \hat \rho^2 ) \beta_2^0 / \lambda$.
As  $| z_2^* | \le 1$ and $\beta_1^* = \beta_1^0- \Delta_1 \ge 0 $, $\beta_2^* =0$,
we see that indeed $\beta^*$ is the solution of the KKT conditions. 

$\circ$ ${\rm Case \ 3 } $:
\fbox {$\lambda/ \hat \varphi^2   > \beta_2^0 + ( \beta_1^0 - \beta_2^0 )/\hat \varphi^2  $}.
With $\Delta_1= \beta_1^0 $ and $\Delta_2= \beta_2^0$ we
obtain
$$ \hat \Sigma \Delta = \begin{pmatrix} \beta_1^0 - \hat \rho \beta_2^0 \cr
- \hat \rho \beta_1^0+ \beta_2^0 \cr \end{pmatrix} = \lambda \begin{pmatrix} z_1^* \cr z_2^* \end{pmatrix} $$
where
$0 <  z_1^* = ( \beta_1^0 - \beta_2^0 + \hat \varphi^2 \beta_2^0 )/\lambda \le 1 $ and
$\lambda z_2^* = \le (- \hat \rho \beta_2^0 + \beta_2^0 )/\lambda = \hat \varphi^2 \beta_2^0 /\lambda \le \lambda$
and $ \lambda z_2^* = (\hat \varphi^2 \beta_1^0 - (\beta_1^0 - \beta_2^0))/\lambda  \ge
-(\beta_1^0 - \beta_2^0 )/\lambda  \ge - 1 + \hat \varphi^2 \beta_2^0/\lambda   \ge - \lambda $. 

Hence the KKT conditions hold for $\beta_1^* = \beta_2^* = 0$. 
\hfill $\sqcap \mkern -12mu \sqcup$ 

\subsection{Proofs for Section \ref{p=s0=2N.section}}

 {\bf Proof of Lemma \ref{eigenvalue.lemma}.} The expression for the minimal eigenvalue
$ \Lambda_{\rm min} (\hat \Sigma)$ is trivial. Then, 
 by orthogonality
 $$  \| X \beta \|_2^2 =\sum_{k=1}^N \| X \beta_{2k-1} + X \beta_{2k} \|_2^2 $$
 and by the arguments of Lemma \ref{hatphi1.lemma} for all $k$
 $$\min_{|\beta_{2k-1} | + |\beta_{2k} | =1 } \| X \beta_{2k-1} + X \beta_{2k} \|_2^2 = \hat \varphi_k^2/2 . $$
 For any vector $v \in \R^N$
 $$ \| v \|_1^2 = \biggl ( \sum_{k=1}^N { | v_k| \hat \varphi_k \over \hat \varphi_k } \biggr ) \le
 \biggl ( \sum_{k=1}^N v_k^2 \hat \varphi_k^2  \biggr ) \biggl ( \sum_{k=1}^N 1/ \hat \varphi_k^2  \biggr )  $$
 and this gives
 $$ \min_{\| v \|_1 = 1} \sum_{k=1}^N v_k^2 \hat \varphi_k^2 = \biggl ( \sum_{k=1}^N 1/ \hat \varphi_k^2
  \biggr )^{-1} =  \| 1 / \hat \varphi^2 \|_1^{-1} . $$
 
 So
 \begin{eqnarray*}
\min_{\| \beta \|_1 =1 }  \| X \beta \|_2^2
 & = & \min_{\| v \|_1 =1} \sum_{k=1}^N \min_{ |\beta_{2k-1} |+ |\beta_{2k} | = v_k } 
 \| X \beta_{2k-1} + X \beta_{2k} \|_2^2\\
 &=& \min_{\| v \|_1 =1 } \sum_{k=1}^N v_k^2 \hat \varphi_k^2/2  =  \| 1 / \hat \varphi \|_1^{-1} /2. 
\end{eqnarray*}

The expression for $\hat \phi^2 ({\cal S})$ follows by similar arguments.
\hfill $\sqcup \mkern -12mu \sqcap$

 {\bf Proof of Lemma \ref{UII-IIIbis.lemma}.} 
 By Lemma \ref{eigenvalue.lemma} the compatibility constant is
 $$ \hat \phi^2 (S_0) = N \| 1/ \hat \varphi^2  \|_1^{-1} . $$
 
This gives  by Lemma \ref{upperbound1.lemma} (recall $|S_0| = 2N$)
$$ \| X (\beta^* - \beta^0 ) \|_2^2 \le 2 \lambda^2 \|  1/ \hat \varphi^2  \|_1 . $$

 On the other hand, by the orthogonality  and the decomposability of the $\ell_1$-norm,
 the Lasso problem can also be decomposed,
 giving in view of Lemma \ref{p=2-exact.lemma}, for each $k$, 
 $$ \| X_{2k-1} (\beta_{2k-1}^* - \beta_{2k-1}^0 ) + X_{2k} (\beta_{2k}^* - \beta_{2k}^0 ) \|_2^2 =
 {2  \lambda^2 \over \hat \varphi_k^2 } $$
 and
 \begin{eqnarray*}
  & \ &  \| X_{2k-1} (\beta_{2k-1}^* - \beta_{2k-1}^0 ) + X_{2k} (\beta_{2k}^* - \beta_{2k}^0 ) \|_2^2  \\
  &+& 
 \lambda (| \beta_{2k-1}^* | +  | \beta_{2k}^* |)\\ &=&
 \lambda (| \beta_{2k-1}^0 | +  | \beta_{2k}^0 |) 
 \end{eqnarray*} 
 where we
 used the assumption $\lambda/ \hat \varphi_k^2 \le \beta_{2k}^0 \le \beta_{2k-1}^0 $. 
 Thus
 $$ \| X ( \beta^* - \beta^0 ) \|_2^2 = 2 \lambda^2 \sum_{k=1}^N 1/ \hat \varphi_k^{2} =
 2 \lambda^2 \|1/  \hat \varphi^2 \|_1  $$
 and
  $$ \| X (\beta^* - \beta^0) \|_2^2 + \lambda \| \beta^* \|_1 = \lambda \| \beta^0 \|_1  . $$
   \hfill $\sqcup \mkern -12mu \sqcap$ 
   
  \subsection{Proof of the lemma in Section \ref{trivial.section}}
   
   {\bf Proof of Lemma \ref{trivial.lemma}.}
   Obviously for all $\beta_{S_0}$
   $$ \arg \min \biggl \{ \| X \beta_{S_0} - X \beta_{-S_0} \|_2^2 : \ \| \beta_{-S_0} \|_1\le 1 \biggr \} = 0 .$$
   So the result of the lemma  follows immediately from Lemmas \ref{hatphi1.lemma} and \ref{p=2-exact.lemma}. 
   \hfill $\sqcup \mkern -12mu \sqcap$
   
   \subsection{Proof of the lemma in Section \ref{s0=2-m0=1.section}}
   
   {\bf Proof of Lemma \ref{s0=2-m0=1.lemma}.} It holds by symmetry arguments that
for all $\beta_3 \in \R$
$$ \min_{| \beta_1| + | \beta_2 | =1 } 
\| X_1 \beta_1 + X_2 \beta_2 - X_3 \beta_3 \|_2^2 =
\| (X_1 +X_2 )/2 - X_3 \beta_3 \|_2^2 . $$
Moreover
\begin{eqnarray*}
\gamma_3& := & \argmin_{\beta_3 \in \R} \| (X_1 +X_2 )/2 - X_3 \beta_3 \|_2^2 \\
& = & \argmin_{\beta_3 \in \R} \| X_1 + X_2 - 2 X_3 \beta_3 \|_2^2 \\
&=& \argmin_{\beta_3 \in \R}   \| (X_1 + X_2) (1-  C \beta_3) - 2 U \beta_3 \|_2^2  \\
& = & \argmin_{\beta_3 \in \R} \biggl \{ 2 \hat \varphi^2 (1-  C \beta_3)^2 + 4 \tau^2 \beta_3^2 \biggr \} \\
& = & { C ( 2 \hat \varphi^2 ) \over 4 \hat \tau^2 + C^2 (2 \hat \varphi^2 ) } . 
\end{eqnarray*}
Since $| \gamma_3 | < 1 $, we conclude that
\begin{eqnarray*}
\hat \phi^2 (S_0) / s_0  &= & \| (X_1 + X_2) / 2 - X_3 \gamma_3 \|_2^2/ n \\
&=& {1 \over 4} \times {  ( 2 \hat \varphi^2 ) ( 4 \hat \tau^2 ) \over
4 \hat \tau^2 + C^2 ( 2 \hat \varphi^2 ) } \\
&=&  { 1 \over 4} \times {  (2 \hat \varphi^2 ) \hat \tau^2 \over \hat \tau^2 + C^2  \hat \varphi^2 /2 } \\
&=&    \hat \varphi^2 \hat \tau^2  /2
\end{eqnarray*} 
where in the last step we used that $\hat \tau^2 + C^2\hat \varphi^2 /2= 1 $. 
Since $s_0=2$ we conclude that $\hat \phi^2 (S_0) = \hat \varphi^2 \hat \tau^2 $. 

To arrive at the second result, we 
write $\beta_1^*= \beta_1^0 - \Delta_1 $ and
$\beta_2^* = \beta_2^0 - \Delta_2 $. 
We have 
$$ \hat \Sigma = \begin{pmatrix}1 & - \hat \rho & C \hat \varphi^2/2 \cr
-\hat \rho & 1 & C \hat \varphi^2/2 \cr C \hat \varphi^2 /2& C \hat \varphi^2 /2& 1 \cr  \end{pmatrix}. $$
For $\beta_3^* = \lambda (2C-1) / \hat \tau^2 $, $\Delta_1= \Delta_2= C \beta_3^* + \lambda / \hat \varphi^2 $
we find
$$ \hat \Sigma  \begin{pmatrix} \Delta_1 \cr \Delta_2 \cr - \beta_3^* \cr \end{pmatrix} =
 \begin{pmatrix} \hat \varphi^2 \Delta_1 - C \hat \varphi^2 \beta_3^* /2\cr
 \hat \varphi^2 \Delta_1 - C \hat \varphi^2 \beta_3^*/2 \cr 
 C \hat \varphi^2 \Delta_1 - \beta_3^* \cr   \end{pmatrix} = \lambda  \begin{pmatrix} 1 \cr 1 \cr 1 \cr \end{pmatrix} 
 .$$
  Since $0 \le \beta_1^* \le \beta_1^0$ and $0 \le \beta_2^* \le \beta_2^0$ and
  $\beta_3^* >0$, the vector $\beta^*$ is indeed the  solution of the KKT conditions. 
  \hfill $\sqcup \mkern -12mu \sqcap$
  
  \subsection{Proof of the lemma in Section \ref{s0=2N-m0=1.section}}
  
  {\bf Proof of Lemma \ref{s0=2N-m0=1.lemma}.}
Let
$$ \gamma:= \arg \min\biggl \{ \| X \beta_{S_0} - X_{2N+1} \beta_{2N+1} \|_2^2   : \ 
\| \beta_{-S_0} \|_1 = 1 , \ | \beta_{2N+1 } | \le 1 \biggr \} . $$
By straightforward computations in a spirit similar to the one used in 
the proof of Lemma \ref{s0=2-m0=1.lemma}, one finds for
$k=1 , \ldots , N$
$$ \gamma_{2k-1} = \gamma_{2k}  =  {  C^2 /(N \| 1/ \hat \varphi^2 \|_1) +
2 \hat \tau^2 \hat \varphi_k^2 /  \| 1/ \hat \varphi^2 \|_1 \over
2C^2 / \| 1/ \hat \varphi^2 \|_1 + 4 \hat \tau^2 } $$ 
and moreover
$$ \gamma_{2N+1} = { 2 C \| 1/ \hat \varphi^2 \|_1 \over
2 C^2\| 1/ \hat \varphi^2 \|_1 + 4 \hat \tau^2 } . $$
Inserting these values one sees
\begin{eqnarray*}
 \hat \phi^2 (S_0) / s_0 &=& \| X \gamma_{S_0} - X_{2N+1} \gamma_{2N+1} \|_2^2 \\
 &=& {  \hat \tau^2 /  \| 1 / \hat \varphi^2 \|_1 \over
 2 \hat \tau^2 + C^2 / \| 1/ \hat \varphi^2 \|_1 } .
 \end{eqnarray*} 
 The second result of the lemma also follows from similar arguments as
 used in the proof of Lemma \ref{s0=2-m0=1.lemma}. The minimizing values
 are
 $$ \beta_{2N+1}^* = \lambda (C -1) / \hat \tau^2 $$
 and for $k=1 , \ldots , N$
 $$\beta_{2k-1}^0 - \beta_{2k-1}^* - C \beta_{2N+1}^*/2 =
 \beta_{2k}^0 - \beta_{2k}^* - C \beta_{2N+1}^*/2 = \lambda / \hat \varphi_k^2 . $$
 \hfill $\sqcup \mkern -12mu \sqcap$
 
 \subsection{Proof of the lemma in Section \ref{s0=2-m0.section}}
 
 {\bf Proof of Lemma \ref{s0=2-m0.lemma}.} This follows by similar arguments as
used in the proof of Lemma \ref{s0=2-m0=1.lemma}.
\hfill $\sqcup \mkern -12mu \sqcap$

\subsection{Proofs for Section \ref{s0=2=m0.section}}

  {\bf Proof of Lemma \ref{goodcomp.lemma}.} 
   We minimize
   $$ \| X_1 \beta_1 + X_2 \beta_2 - X_3 \beta_3 - X_4 \beta_4 \|_2^2 $$
   over $|\beta_1 |+ |\beta_2 |=1$ and $|\beta_3 |+ |\beta_4 |\le 1$. 
   It holds that
   \begin{eqnarray*}
 & \ &    \| X_1 \beta_1 + X_2 \beta_2 - X_3 \beta_3 - X_4 \beta_4 \|_2^2 \\
 &= & \| X_1 ( \beta_1 - C(\beta_3 + \beta_4)/2) + X_2 ( \beta_2- C(\beta_3 + \beta_4)/2) \|_2^2\\
 & +&
 (\beta_3 -\beta_4)^2  (1-    C^2 \hat \varphi^2 /2- \hat \tau^2  ) + (\beta_3+ \beta_4)^2\hat \tau^2 
    \end{eqnarray*}
   This implies $\beta_3=\beta_4$.
   So we minimize
  $$  \| X_1 ( \beta_1 - C\beta_3 ) + X_2 ( \beta_2- C\beta_3 ) \|_2^2 + 
  (2\beta_3)^2 \hat \tau^2  . $$
  By symmetry arguments, we know $\beta_1 = \beta_2 $ say both $+1/2$.
  Then we need to minimize
  $$  \| (X_1+X_2 )( 1/2 - C\beta_3 )  \|_2^2 + 
  (2\beta_3)^2 \hat \tau^2   = (1/2 - C \beta_3)^2 2\hat \varphi^2  + (2\beta_3)^2 \hat \tau^2. $$
  The minimizing value for $2\beta_3$ is 
  $$ 2 \gamma_3 = {1 \over 2} { 2 C (2  \hat \varphi^2 )\over
 C^2 ( 2 \hat \varphi^2 ) + 4 \hat \tau^2 } . $$
  In other words
  $$ \min \biggl \{ \| X \beta \|_2^2 : \ \| \beta_{S_0} \|_1=1 , \ \| \beta_{-S_0} \|_1 \le 1 \biggr \} =
  {1 \over 4} { 2 \hat \varphi^2 \hat \tau^2 \over C^2 (2 \hat \varphi^2) + \hat \tau^2 } . $$
  Hence
  $$ \hat \phi^2 (S_0) = {1 \over 2}   { 2 \hat \varphi^2 \hat \tau^2 \over C^2 (2 \hat \varphi^2) + \hat \tau^2 }. $$
  
 For the second result, we check the KKT conditions with
 $\beta_3^*= \beta_4^* $ and $\Delta_1 - C \beta_3^*= \Delta_2-C \beta_3^* = \lambda
   / \hat \varphi^2 $ where $\Delta_1= \beta_1^0 - \beta_1^*$ and $\Delta_2= \beta_2^0 - \beta_2^* $. 
     It holds that
   $$ \hat \Sigma= \begin{pmatrix} 1 & - \hat \rho & C\hat \varphi^2/2 & C\hat\varphi^2/2 \cr 
   - \hat \rho & 1 & C \hat \varphi^2 /2& C \hat\varphi^2/2 \cr 
  C \hat \varphi^2 /2& C \hat\varphi^2 /2& 1 &   C^2  \hat \varphi^2 + 2 \hat \tau^2 -1\cr 
    C \hat \varphi^2 /2& C \hat\varphi^2 /2&    C^2 \hat \varphi^2 + 2 \hat \tau^2 -1 &1 \cr 
   \end{pmatrix} .$$
   Hence
   \begin{eqnarray*}
    \hat \Sigma \begin{pmatrix} \Delta_1 \cr  \Delta_2 \cr- \beta_3 \cr -\beta_4 \cr 
   \end{pmatrix} &=& \begin{pmatrix}  C\hat \varphi^2 \Delta_1/2  -
  C \hat \varphi^2 \beta_3 \cr   C\hat \varphi^2 \Delta_1/2  -
   C \hat \varphi^2 \beta_3 \cr  
   C \hat \varphi^2 \Delta_1  - ( C^2 \hat \varphi^2 + 2\hat \tau^2) \beta_3 \cr
    C \hat \varphi^2 \Delta_1  - ( C^2  \hat \varphi^2 + 2\hat \tau^2 ) \beta_3 \cr
   \end{pmatrix}\\
   & = & \begin{pmatrix}  C\hat \varphi^2 ( \lambda/ \hat \varphi^2 + C \beta_3 )/2 -
  C \hat \varphi^2 \beta_3 \cr   C\hat \varphi^2 ( \lambda/ \hat \varphi^2 +  C \beta_3 ) /2-
  C \hat \varphi^2 \beta_3 \cr  
   C \hat \varphi^2 ( \lambda/ \hat \varphi^2 +  C\beta_3) -( C^2\hat \varphi^2 + 2\hat \tau^2 ) \beta_3 \cr
    C \hat \varphi^2 ( \lambda/ \hat \varphi^2 + C \beta_3 ) - ( C^2\hat \varphi^2 + 2\hat \tau^2  ) \beta_3 \cr
   \end{pmatrix}\\
    & = & \begin{pmatrix}   \lambda 
   \cr   \lambda \cr  
    C \lambda - 2\hat \tau^2 \beta_3 \cr
    C \lambda  - 2\hat \tau^2 \beta_3 \cr
   \end{pmatrix}  = \lambda \begin{pmatrix}1 \cr  1 \cr  1 \cr  1 \cr \end{pmatrix} .  
   \end{eqnarray*} 
  Thus $\beta^*$ is the solution of the KKT conditions. 
    \hfill $\sqcup \mkern -12mu \sqcap$
    
     {\bf Proof of Lemma \ref{goodlasso2.lemma}.}
     Along similar lines as in the proof of  Lemma \ref{goodcomp.lemma}, one
     finds $\beta_3^*= \beta_4^*$ and 
     $$  \beta_2^0 - \beta_2^* - C \beta_3^*  = \beta_1^0 - \beta_1^* - C \beta_3^* =
     \lambda / \hat \varphi^2  $$
     but now $\beta_3^*$ is the largest possible value such that
     $  \beta_2^0 - \beta_2^* \le  \beta_2^0 $. It follows
     that $\beta_1^* = \beta_1^0 -
     \beta_2^0 $, $ \beta_2^* = 0$, $ \beta_3^* = ( \beta_2^0 - \lambda / \hat \varphi^2 ) / C $.
     \hfill $\sqcup \mkern -12mu \sqcap$
     
      {\bf Proof of Lemma \ref{goodlasso3.lemma}.}
    The Gram matrix is now
    $$ \hat \Sigma = \begin{pmatrix}1 & - \hat \rho & \hat \varphi^2 / 2 & \hat \varphi^2 / 2 \cr 
    - \hat \rho & 1 &  \hat \varphi^2 / 2 & \hat \varphi^2 / 2 \cr 
    \hat \varphi^2 / 2 & \hat \varphi^2 / 2 & 1 & \hat \varphi^2 -1 \cr
    \hat \varphi^2 / 2 & \hat \varphi^2 / 2 & \hat \varphi^2 -1 & 1 \cr 
    \end{pmatrix} .$$
    Hence, with $\Delta_1 = \Delta_2= \beta_3^* + \lambda / \hat \varphi^2 $ we find
    \begin{eqnarray*}
     \hat \Sigma \begin{pmatrix} \Delta_1 \cr \Delta_2 \cr - \beta_3^* \cr - \beta_3^* \cr \end{pmatrix} &=&
    \begin{pmatrix} \hat \varphi^2 ( \beta_3^* + \lambda / \hat \varphi^2 ) - 
   \hat \varphi^2 \beta_3^* \cr
   \hat \varphi^2 ( \beta_3^* + \lambda / \hat \varphi^2 ) - 
   \hat \varphi^2 \beta_3^* \cr  \hat \varphi^2 ( \beta_3^* + \lambda / \varphi^2 ) -
   \beta_3^* - ( \hat \varphi^2 -1) \beta_3^* \cr 
   \hat \varphi^2 ( \beta_3^* + \lambda / \varphi^2 ) -
   \beta_3^* - ( \hat \varphi^2 -1) \beta_3^* \cr \end{pmatrix} \\
   &=& \lambda \begin{pmatrix} 1 \cr 1 \cr 1 \cr 1 \cr \end{pmatrix} . 
   \end{eqnarray*} 
   Since for $ 0 \le \beta_3^* \le \beta_2^0 - \lambda / \hat \varphi^2 $ it holds 
   that
   $\beta_1^* = \beta_1^0 - \beta_3^* - \lambda /\hat  \varphi^2 \ge 0 $, 
   $\beta_2^* = \beta_2^0 - \beta_3^* - \lambda /\hat  \varphi^2 \ge 0$ and $\beta_3^* \ge 0$,
   the vector $\beta^*$ is indeed the solution of the KKT conditions.
   With this value one finds the result for the prediction error and the bound for
   $\| \beta_{-S_0}^* \|_1 $.
   \hfill $\sqcup \mkern -12mu \sqcap$
   
   \subsection{Proof of the lemma in Section \ref{s0=m0=2N.section}}
   
  {\bf Proof of Lemma \ref{s0=m0=2N.lemma}.} This follows from the same
  arguments as used in the proofs of Lemmas \ref{eigenvalue.lemma} and  \ref{goodcomp.lemma}.
  \hfill $\sqcup \mkern -12mu \sqcap$
  
  \subsection{Proofs for Section \ref{Furthers0=2.section}}
   
   {\bf Proof of Lemma \ref{m=s0=2.lemma}.}
Observe first that 
\begin{eqnarray*}
 \hat \rho &=&1-  2 \gamma_{-S_0}^T \hat \Sigma_{-S_0, -S_0} \gamma_{-S_0}   \\
 &=& 1- 2 \biggl ( \gamma_3^2 + (1- \gamma_3)^2 - 2 \gamma_3 (1- \gamma_3) \hat \theta \biggr )  \\
 &=& 1- 2 \biggl ( 1 - 2 \gamma_3 (1- \gamma_3) (1+ \hat \theta)  \biggr )\\
 &=& = -1 +  4 \gamma_3 (1- \gamma_3) (1+ \hat \theta) .
 \end{eqnarray*}
Therefore
\begin{eqnarray*}
\hat \varphi^2 & = & 1
- \hat \rho\\ & =&  2- 4 \gamma_3 (1- \gamma_3) (1+ \hat \theta) \\
&=&
2(1- 4 \gamma_3 (1- \gamma_3)) + 4 \gamma_3 (1- \gamma_3)\hat \psi^2 . 
 \end{eqnarray*}
 So $\hat \varphi^2 >  \hat \psi^2 $.
 
 The Gram matrix is now
 $$ \hat \Sigma = \begin{pmatrix} 1  & - \hat \rho & \hat \varrho_3&
\hat \varrho_4   \cr 
 -\hat \rho & 1 & \hat \varrho_3&
\hat \varrho_4 \cr  
\hat \varrho_3 & \hat \varrho_3    &1 & -\hat \theta \cr
\hat \varrho_4 & \hat \varrho_4    & -\hat \theta & 1  \cr \end{pmatrix} $$
  where $\hat \varrho_3 :=
  \gamma_3- (1- \gamma_3) \hat \theta $ and $\hat \varrho_4= (1- \gamma_3) - \gamma_3 \hat \theta$.
  Then for $\beta_1^* = \beta_1^0 - \Delta_1:= \beta_1^0 - \beta_2^0$, $\beta_2^* = 
  \beta_2 - \Delta_2 :=0$, $\beta_3^*= 2 \gamma_3
  \beta_2^0- \lambda/ \hat \psi^2 $ and $\beta_4^* = 2 (1- \gamma_3) \beta_2^0 -
  \lambda / \hat  \psi^2 $ we get
  $$ \hat \Sigma \begin{pmatrix}\Delta_1 \cr \Delta_2 \cr - \beta_3^* \cr - \beta_4^* \cr  \end{pmatrix} =
  \begin{pmatrix} \hat \varphi^2 \beta_2^0 - \hat \varrho_ 3\beta_3^*+ \hat\varrho_4 \beta_4^*  \cr 
  \hat \varphi^2 \beta_2^0 - \hat \varrho_3 \beta_3^* +\hat \varrho_4 \beta_4^* \cr 
  2  \beta_2^0 \hat \varrho_3   + \beta_3^* - \hat \theta \beta_4^* \cr 
    2 \beta_2^0 \hat \varrho_4- \hat \theta  \beta_3^*  + \beta_4^* \cr 
  \end{pmatrix}= \lambda \begin{pmatrix} 1 \cr 1 \cr 1 \cr 1 \cr \end{pmatrix} . 
  $$
  So $\beta^*$ is the solution of the KKT conditions.
  \hfill $\sqcup \mkern -12mu \sqcap$

   {\bf Proof of Lemma \ref{phi=rho.lemma}.}  
 It is straightforward to calculate
  $$\hat \varphi^2  = C^2 \hat \psi^2. $$
  To find $\hat \phi^2 (S_0)$ we minimize
  $$ \| X_1\beta_1 +  X_2 \beta_2 - X_{-S_0} \beta_{-S_0} \|_2^2 $$
  over $0 < \beta_1 = 1 \beta_2 < 1$ and $\| \beta_{-S_0} \|_1 \le 1 $.
  Symmetry arguments yield $\beta_1=1/2$.
  We then minimize
  $$ \|   X_3 + X_4  -  X_{S_0} \beta_{-S_0}) \|_2^2 $$ over
  $\| \beta_{-S_0} \|_1 \le 1$. This gives that the entries
  in $\beta_{-S_0}$ are equal to $1/2$ and hence $\hat \phi^2 (S_0) = (C-1)^2\hat  \psi^2 $.
  
   In view of Lemmas \ref{p=2-exact.lemma} and \ref{phi=rho.lemma},
   it suffices to show that $\beta_{-S_0}^* = 0 $ corresponds to the unique solution of
   the KKT conditions.
   We have with $\Delta_1=\Delta_2= \lambda/ \hat \varphi^2$
   \begin{eqnarray*}
     \hat \Sigma_{-S_0, S_0} \begin{pmatrix}\Delta_1 \cr \Delta_2 \cr 
   \end{pmatrix}&=&   \begin{pmatrix}  C \hat \psi^2 / 2  & C \hat \psi^2 / 2  \cr C \hat \psi^2 / 2 & C \hat \psi^2 / 2 \cr
    \end{pmatrix} \begin{pmatrix}\lambda/ \hat \varphi^2 \cr  \lambda/ \hat \varphi^2 \cr \end{pmatrix} \\
   &=& \begin{pmatrix} {C\lambda \hat \psi^2 /  \hat \varphi^2  } \cr { C \lambda \hat \psi^2 / \hat \varphi^2  } \cr \end{pmatrix} \\ &  =&  \lambda \begin{pmatrix} 1/C \cr 1/C \cr \end{pmatrix}, 
   \end{eqnarray*}
  since $\hat \varphi^2  =C^2 \hat \psi^2  $. So the KKT conditions are satisfied, with
  $z_{-S_0}^* = 1 / C $. 
  The solution is unique because for $\gamma_{-S_0} = (C  , C )^T/2 $ it holds that $\| \gamma_{-S_0} \|_1 > 1 $. 
    \hfill $\sqcup \mkern -12mu \sqcap $
    
    {\bf Proof of Lemma \ref{2comp>0.lemma}.} 
  Note first that indeed $\hat \rho =1-  2 C^2  \| \gamma_{-S_0}  \|_2^2 >0$ since $2 C^2 \| \gamma_{-S_0} \|_2^2 < 1$.
  It follows that $\hat \varphi^2 = 2 C^2  \| \gamma_{-S_0} \|_2^2 $.
  We have
  $$ \| \beta_1 X_1 + \beta_2 X_2 - X_{-S_0} \beta_{-S_0} \|_2^2 $$
  is minimized over $|\beta_1 | + |\beta_2 | =1$ at $\beta_1=\beta_2 = 1/2$ and
  $$ \| X_{-S_0} (C \gamma_{-S_0}  - \beta_{-S_0})  \|_2^2 = \| C \gamma_{-S_0}  - \beta_{-S_0}  \|_2^2
 . $$
  To obtain the prediction error, in view of Lemmas \ref{p=2-exact.lemma} and \ref{phi=rho.lemma},
  it suffices to show that $\beta_{-S_0}^* =0 $ and $\Delta_1= \beta_1^0- \beta_1^*=
  \lambda / \hat \varphi^2 $, $\Delta_2= \beta_1^0- \beta_1^*=\lambda / \hat \varphi^2 $
  is the unique solution of the KKT conditions.
  We have 
  $$ \hat \Sigma_{-S_0,  S_0} \begin{pmatrix} \Delta_1 \cr \Delta_2\cr  \end{pmatrix}=    
 C \gamma_{-S_0}  (\Delta_1 + \Delta_2 ) = 2 C \lambda \gamma_{-S_0}  / \hat \varphi^2 = \lambda \gamma_{-S_0} / C \| \gamma_{-S_0}  \|_2^2 =
  \lambda z_{-S_0}^* ,  $$
  where $\| z_{-S_0}^* \|_{\infty}  = \| \gamma_{-S_0} \|_{\infty}  /  (C\| \gamma_{-S_0} \|_2^2 ) \le 
  1 $.   The solution is unique because $\| \gamma_{-S_0} \|_1 > 1 $.
  Another way to see it is by noting that for any  $\| z_{S_0} \|_{\infty} \le 1 $
  $$ \| \hat \Sigma_{-S_0, S_0} \hat \Sigma_{S_0 , S_0}^{-1} z_{S_0} \|_{\infty} =
  \| \gamma_{-S_0} (z_1 + z_2 )  \|_{\infty} / \hat \varphi^2 = \| \gamma_{-S_0 } \|_{\infty} / (C\| \gamma_{-S_0} \|_2^2 ) < 1 $$
  i.e., the irrepresentable condition holds.

\bibliographystyle{plainnat}
\bibliography{reference}
\end{document}